\documentclass[
  draft,
leqno,   a4paper,
]{article}
\usepackage{ makeidx,latexsym, color, amssymb,  amsmath,stmaryrd}

\usepackage{amsfonts, amsmath, amssymb, amsthm, amsxtra}
\usepackage{graphics}

\theoremstyle{plain}
\newtheorem{Theorem}{Theorem}[section]
\newtheorem{Prop}[Theorem]{Proposition}
\newtheorem{Lemma}[Theorem]{Lemma}
\newtheorem{Corollary}[Theorem]{Corollary}

\theoremstyle{definition}

\newtheorem{Remark}[Theorem]{Remark}

\newtheorem{Example}[Theorem]{Example}

\begin{document}
\title{Meromorphic Continuation of Dynamical Zeta Functions via Transfer Operators}
\author{J. Hilgert\thanks{Universit\"{a}t Paderborn, Institut f\"{u}r Mathematik, Warburger Str. 100, D-33098 Paderborn, Germany; hilgert@math.upb.de}{ }\footnote{corresponding author.} \and   F. Rilke\thanks{Universit\"{a}t Paderborn, Institut f\"{u}r Mathematik, Warburger Str. 100, D-33098 Paderborn, Germany; rilke@math.upb.de}{ }\footnote{Supported by the Deutsche Forschungsmeinschaft
through  the research projects ``Gitterspinsysteme'' and
``Gitterspinsysteme II''.} }
\date{}
\maketitle
%
%
\def\id{\mathop{\text{\rm id}}}
\renewcommand{\Im}[1]{\ensuremath  \mathbf{I\!m}{(#1)}}
\renewcommand{\Re}[1]{\ensuremath  \mathbf{R\!e}{(#1)}}
\newcommand{\Abl}[2]{\frac{d{#1}}{d{#2}}}
\newcommand{\ZAbl}[2]{\frac{d^2{#1}}{d{#2}^2}}
\newcommand{\PAbl}[2]{\frac{\partial{#1}}{\partial{#2}}}
\newcommand{\ZPAbl}[3]{\frac{\partial^2{#1}}{\partial{#2}\partial{#3}}}
\newcommand{\ZPAbla}[2]{\frac{\partial^2{#1}}{\partial{#2}^2}}

 \newcommand{\euclsp}[2]{\left(\!\!\begin{array}{c|c}\!{#1}\!&\!{#2}\!\end{array}\!\!\right)}
 \newcommand{\hermsp}[2]{\left<\!\!\begin{array}{c|c}\!{#1}\!&\!{#2}\!\end{array}\!\!\right>}

\def\ssarr{\hbox to 30pt{\rightarrowfill}}
\def\sarr{\hbox to 40pt{\rightarrowfill}}
\def\arr{\hbox to 55pt{\rightarrowfill}}
\def\larr{\hbox to 55pt{\leftarrowfill}}
\def\Arr{\hbox to 80pt{\rightarrowfill}}
\def\doublemapright#1{\smash{\mathop{\doublearr}\limits^{#1}}}
\def\mapdown#1{\Big\downarrow\rlap{$\vcenter{\hbox{$\scriptstyle#1$}}$}}
\def\lmapdown#1{\llap{$\vcenter{\hbox{$\scriptstyle#1$}}$}\Big\downarrow}
\def\mapright#1{\smash{\mathop{\arr}\limits^{#1}}}
\def\mapleft#1{\smash{\mathop{\larr}\limits^{#1}}}
\def\ssmapright#1{\smash{\mathop{\ssarr}\limits^{#1}}}
\def\smapright#1{\smash{\mathop{\sarr}\limits_{#1}}}
{}\def\lmapright#1{\smash{\mathop{\arr}\limits_{#1}}}
\def\Mapright#1{\smash{\mathop{\Arr}\limits^{#1}}}
\def\doublearr{\hbox to 110pt{\rightarrowfill}}
\def\lMapright#1{\smash{\mathop{\Arr}\limits_{#1}}}
\def\mapup#1{\Big\uparrow\rlap{$\vcenter{\hbox{$\scriptstyle#1$}}$}}
\def\lmapup#1{\llap{$\vcenter{\hbox{$\scriptstyle#1$}}$}\Big\uparrow}
\def\vline{\hbox{\Bigg\vert}}
%
\newcommand{\prtbox}{ \protect\phantom{\fbox{1}} }
\newcommand{\prtboxa}{ \protect{\fbox{1}} }
\newcommand{\prtboxb}{ \protect{\fbox{2}} }
\newcommand{\prtboxc}{ \protect{\fbox{3}} }
\newcommand{\xyvektor}{\left(\!\begin{array}{c} x\\ y\end{array}\!\right)}
\newcommand{\UP}[1]{\textup{#1}}
\renewcommand{\emph}[1]{\textbf{#1}}
\newcommand{\Acal}{\mathcal{A}} 
\newcommand{\Bcal}{\mathcal{B}} 
\newcommand{\Ccal}{\mathcal{C}} 
\newcommand{\Dcal}{\mathcal{D}} 
\newcommand{\Ecal}{\mathcal{E}} 
\newcommand{\Fcal}{\mathcal{F}} 
\newcommand{\Gcal}{\mathcal{G}} 
\newcommand{\Hcal}{\mathcal{H}} 
\newcommand{\Kcal}{\mathcal{K}} 
\newcommand{\Lcal}{\mathcal{L}} 
\newcommand{\Mcal}{\mathcal{M}} 
\newcommand{\Ocal}{\mathcal{O}} 
\newcommand{\Pcal}{\mathcal{P}} 
\newcommand{\Rcal}{\mathcal{R}} 
\newcommand{\Scal}{\mathcal{S}} 
\newcommand{\Ycal}{\mathcal{Y}} 

\newcommand{\A}{\mathbb{A}}
\newcommand{\B}{\mathbb{B}}
\newcommand{\C}{\mathbb{C}}    
\newcommand{\D}{\mathbb{D}}
\newcommand{\E}{\mathbb{E}}
\renewcommand{\H}{\mathbb{H}}    
\newcommand{\K}{\mathbb{K}}
\newcommand{\LL}{\mathbb{L}}
\newcommand{\MM}{\mathbb{M}}
\newcommand{\N}{\mathbb{N}} 
\newcommand{\No}{{\mathbb{N}_{0} }} 
\newcommand{\PP}{\mathbb{P}}
\newcommand{\Q}{\mathbb{Q}}
\newcommand{\R}{\mathbb{R}} 
\renewcommand{\S}{\mathbb{S}}    
\newcommand{\T}{\mathbb{T}}
\newcommand{\Z}{\mathbb{Z}} 
\newcommand{\dist}{\mbox{\UP{dist}}}
\newcommand{\spec}{\mbox{\UP{spec}}}
\newcommand{\einhalb}{\frac{1}{2}}
\newcommand{\dreihalb}{\frac{3}{2}}
\newcommand{\einviertel}{\frac{1}{4}}

\newcommand{\dom}{\mbox{\UP{dom }}}
\newcommand{\sign}{\mbox{\UP{sign}}}
\newcommand{\trace}{\mbox{\UP{trace }}}
\newcommand{\Tracec}{\mbox{\UP{Trace}}^c\,}
\newcommand{\Tracef}{\mbox{\UP{Trace}}^\flat\,}
\newcommand{\Traces}{\mbox{\UP{Trace}}^\sharp\,}
\newcommand{\Trace}{\mbox{\UP{Trace }}}
\newcommand{\wurzel}[1]{\sqrt{#1}}
\newcommand{\Mat}{\mbox{\UP{Mat}}}
\newcommand{\Gl}{\mbox{\UP{Gl}}}
\newcommand{\PD}{\mbox{\UP{PD}}}
\newcommand{\Sl}{\mbox{\UP{Sl}}}
\newcommand{\Sym}{\mbox{\UP{Sym}}}
\newcommand{\Sp}{\mbox{\UP{Sp}}}
\newcommand{\Lin}{\mbox{\UP{Lin}}}
\newcommand{\End}{\mbox{\UP{End}}}
\newcommand{\sMat}{\mbox{\small\UP{Mat}}}
\newcommand{\sGl}{\mbox{\small\UP{Gl}}}
\newcommand{\sPD}{\mbox{\small\UP{PD}}}
\newcommand{\sSl}{\mbox{\small\UP{Sl}}}
\newcommand{\sSym}{\mbox{\small\UP{Sym}}}
\newcommand{\sSp}{\mbox{\small\UP{Sp}}}
\newcommand{\sLin}{\mbox{\small\UP{Lin}}}
\newcommand{\sEnd}{\mbox{\small\UP{End}}}
\newcommand{\arccosh}{\mbox{\UP{arccosh}}}
\newcommand{\Per}{\mathrm{Per}}
\newcommand{\Pol}{\mathrm{Pol}}
\newcommand{\diag}{\mathrm{diag}}
\newcommand{\range}{\mathrm{range}}
\newcommand{\rank}{\mathrm{rank}}
\newcommand{\co}{\mathrm{co}}
\renewcommand{\span}{\mathrm{span}}
\newcommand{\pr}{\mathrm{pr}}
\newcommand{\vol}{\mathrm{vol}}
\newcommand{\conc}{\mathrm{conc}}
\newcommand{\period}{\mathrm{period}}
\newcommand{\opnorm}{{\mbox{\tiny \UP{op}}}}
\newcommand{\rhospec}{\rho_{\mbox{\small\UP{spec}}}}

\newcommand{\Fix}{\mathrm{Fix}}
\newcommand{\ev}{\mathrm{ev}}
\newcommand{\ds}{\displaystyle}
\renewcommand{\theenumi}{\UP{\roman{enumi}}}
\renewcommand{\labelenumi}{\UP{(}\theenumi\UP{)}}
 \newcommand{\zn}{Z_n}
%
%
%
 \begin{abstract}
 We describe a general method to prove meromorphic continuation of dynamical
 zeta functions to the entire complex plane under the condition that the corresponding
 partition functions are given via a dynamical trace formula from a family of transfer operators.
 Further we give general conditions for the partition functions associated with general spin chains to be
 of this type and provide various families of examples for which these conditions are satisfied.

  \end{abstract}
%

\section*{Introduction}

The \textit{dynamical zeta functions} of interest in this paper are generating functions
of the form
\begin{equation}\label{zeta-r-def}
\zeta_R(z):=\exp\!\Big(\sum_{n=1}^\infty\frac{z^n}{n}\, \zn\Big).
\end{equation}
associated with sequences $\big(\zn\big)^{}_{n\in \N}$ in $\C$.
If the
 $\zn$ arise as {\em partition functions} of a dynamical system, the notion of
 dynamical zeta functions or \textit{Ruelle zeta functions} has  been introduced by Ruelle in ~\cite{Ru76} and \cite{Ru76a}.
Naming and special form of these functions are motivated by concrete examples from statistical mechanics.
The partition function  encodes the statistical properties of a system in thermodynamic equilibrium.
It depends on the temperature, the volume, and the microstates of a finite number
of particles.
We   will consider  partition functions
 of the form
\begin{equation}\label{zn-r-def}
\zn
=\sum_{\nu=0}^{n_1} (-1)^\nu \,\trace G_{\nu}^n\quad\text{or}\quad
\zn= \det(1-\Lambda^n)\,\trace G^n
\end{equation}
for large $n$.
In the first case one has a (possibly) infinite   family of compact operators $G_{\nu}$, called
{\em transfer operators}, such that the corresponding Schatten norms satisfy
$\sum_{\nu=0}^{n_1}\| G_{\nu}^{n_o}\|_{\Scal_{1}(\Hcal_\nu)}<\infty$.  In the second  case
one has as a transfer operator $G$ together with an auxiliary operator
$\Lambda$. We will refer to formulae of the type (\ref{zn-r-def}) as {\em dynamical trace formulae}.
Examples of such dynamical zeta functions derived from trace formulae of the
type (\ref{zn-r-def}) have been treated repeatedly in the literature,
see e.g. \cite{KrW41}, \cite{M76}, \cite{V76},  \cite{VM77}, \cite{M80a}, \cite{Mo89},
\cite{M91}, \cite{HM02},   \cite{HM04}, \cite{R03}, \cite{R06}.

Unlike other kinds of zeta functions such as Riemann's, Selberg's, or Artin's zeta function,
our dynamical zeta function is an exponential of a power series, hence itself a power series.
Considering $s\mapsto \zeta_R(e^{-s})$ one obtains a function which is holomorphic in a  right half plane,
provided $\zeta_R$ has a non-zero radius of convergence.
Zeta functions typically occur as a kind of generating functions for collections of objects
like prime numbers or prime geodesics and  it is natural to ask for their analytic properties.
In particular, in order to prove asymptotic results for counting functions it is important to
know whether the zeta functions have meromorphic continuations to a larger set or even to the entire
complex plane. Moreover, one wants to obtain information on the location of the poles. We will show that
this kind of information can indeed be provided {\it if} the partition function can be written via
dynamical trace formulae. More precisely, we show that such dynamical zeta functions are
quotients of regularized Fredholm determinants, one factor for each transfer operator
(see Theorem~\ref{main-thm-1} and Theorem~\ref{main-thm-3}).
Thus they have representations as Euler products, and it is possible to give a spectral interpretation
for its zeros and poles.

Results of the kind sketched above become interesting only if one has a sufficient supply of examples
of partition functions with the desired properties. We provide a general principle how to
construct such examples from (classical) {\it spin chains}. To describe these, we consider a Hausdorff space  $F$
equipped with a finite  measure $\nu$, and  let
$\A:F\times F\to\{0,1\}$ be a $\nu\otimes \nu$-measurable function, which we call a
{\it transition matrix}. Then
$\Omega_\A:=\{\xi\in F^\N\, |\, \A(\xi_i,\xi_{i+1})=1\ \forall i\}$ will be referred to as
a {\it configuration space}. Further we fix an {\it interaction} $A\in \Ccal_b(\Omega_\A)$.
These data, together with the left shift $\tau:F^\N\to F^\N,\ (\tau \xi)_k:=\xi_{k+1}$,
are called a {\it spin chain} or, more technically, a {\it one-sided one-dimensional matrix subshift}.
With such a matrix subshift
we associate a {\it dynamical partition function}
\begin{equation}\label{partition-fct}
Z_n(A)
:=\int_{F^n}\prod_{i=1}^n
\A_{x_i,x_{i+1 }}\exp\!\Big(\sum_{k=0}^{n-1} A(\tau^k(\overline{x_1\ldots x_n}))\Big)
\, d\nu^n(x_1,\ldots,x_n)
,
\end{equation}
where  $\overline{x_1\ldots x_n}:=(x_1,\ldots,x_n,x_1,\ldots,x_n,\ldots)$ and   $x_{n+1}:=x_{1}$.
Consider the special case of  $A=0$,
and let $\rho_n:F^\N\to F^n, \, \xi\mapsto (\xi_1,\ldots, \xi_n)$
 be the projection. Then
 $$
\zn(0)= \nu^n\big(\rho_n(\{\xi\in \Omega_\A\mid \tau^n\xi=\xi\})\big) ,$$
which
measures the number of closed $\tau$-orbits in $\Omega_\A$ with
period length $n$ with respect to the a priori measure $\nu$. In
particular, if the system is a   {\it full shift}, i.\!~e., $\A
\equiv  1$, then $\zn(0)=\nu(F)^n$. For general non-interacting
matrix subshifts it is possible to show (see Proposition~\ref{non-interacting-case-cor})
that there is an operator
$\Gcal_\A$ such that $\zn(0)=\trace \Gcal_\A^n$. If the interaction
$A$ is non-zero we have to make more assumptions in order to
guarantee such trace representations of the partition function. See Theorem~\ref{thm-1}
for a precise formulation. Its proof depends on a two special types of trace formulae. One
(see Lemma~\ref{lemma-T-HS}) is elementary but subtle and deals with iterates of averages
of Hilbert-Schmidt operators. The other (see Theorem~\ref{spurformel-fockinfty-thm})
deals with composition operators on Fock spaces and is based on a fixed point
formula of Atiyah and Bott.

In order to satisfy the hypotheses of Theorem~\ref{thm-1} one has to
verify certain estimates which boil down to asking for a rather rapid decay of the interactions
(see Theorem~\ref{thm-2}).
Nevertheless, the theorem is strong enough to produce, among others, all the corresponding
results scattered in the literature quoted above (see Example~\ref{classification-example}).

The paper is organized as follows.
In Section~\ref{traces-determinants-sec} we recall some key definitions and
provide a number of technical results concerning traces and determinant in infinite dimensions
used later in the paper.
In  Section~\ref{MeromCont} we show how to continue dynamical zeta functions meromorphically if
 a dynamical trace formula holds.
Section~\ref{fock-sec} contains a description of the composition operators that are instrumental in
the construction of dynamical trace fromulae for spin chains in Section~\ref{dyntraceformulae}.
The underlying trace formulae for these
operators are proven in Section~\ref{traces-sec}.
In Section~\ref{ising-type-sec} we apply our results to Ising type spin chains
and give examples for interactions for which our results can be applied.

\section{Traces and Determinants}
    \label{traces-determinants-sec}

 Given a compact operator $K$ on a Hilbert space $\Hcal$ we will
denote by $(\lambda_j(K))_{j\in \N}$ the sequence of its eigenvalues
counted with multiplicities. Let   $(s_j(K))_{j\in \N}$ be the
sequence of the singular numbers of $K$, i.\!~e., the eigenvalues of
the positive compact operator $|K|=\sqrt{K^\star K}$. For $1\leq
p<\infty$ the \textit{Schatten class} $\Scal_p(\Hcal)$ is defined as
the space of all operators $K$ such that
$$
\|K\|_{\Scal_p(\Hcal)}:=\|(s_n(K))_{n\in\N}\|_{\ell^p(\N)}<\infty.
$$
The Schatten classes $\Scal_p(\Hcal)\subset \End(\Hcal)$ for $1\leq p<\infty$ are
embedded subalgebras,  i.\!~e., they satisfy
    $$\|A\|_{\sEnd(\Hcal)}\leq \|A\|_{\Scal_p(\Hcal)},\quad
    \|AB\|_{\Scal_p(\Hcal)} \leq \|A\|_{\Scal_p(\Hcal)}\,\|B\|_{\Scal_p(\Hcal)},
    $$
and have  the approximation property, i.\!~e.,
 the finite rank operators are dense  with respect to
 $\|\,\cdot\,\|_{\Scal_p(\Hcal)}$ (cf. \cite[Thm. XI.11.1]{GoGoKr00}).
For any
$A_1,\ldots,A_{\lceil p \rceil}\in \Scal_p(\Hcal)
$
one has  (cf. \cite[Thm. IV.11.2]{GoGoKr00})
\begin{equation}
| \trace (A_1\cdots A_{\lceil p \rceil})|
\leq  \|A_1\cdots A_{\lceil p \rceil}\|_{\Scal_1(\Hcal)}\leq
 \prod_{j=1}^{\lceil p \rceil}\|A_j\|_{\Scal_p(\Hcal)}.
\label{est-trace-m}%
\end{equation}
By
\cite[Thm.~XI.1.1]{GoGoKr00}
this estimate    implies that
for any $n_o\in \N_{\geq p}$ the \textit{$n_o$-regularized determinant}
\begin{equation}
    \label{u-det-eq}%
\det\nolimits_{n_o}(1-F):=
\det(1-F)\, \exp\!\Big(
\sum_{k=1}^{n_o-1}\frac{1}{k}\,\trace F^k\Big)
\end{equation}
defined on finite rank operators admits a
 continuous extension to  $\Scal_p(\Hcal)$, also denoted by
 $A\mapsto \det_{n_o}(1-A)$.
The following assertions are true on the level of finite rank operators and
can be extended to $\Scal_p(\Hcal)$ by continuity  (see \cite[Thm.~XI.2.1]{GoGoKr00}).
\begin{Lemma}
    \label{gogokr-lemma-m}%
Let $\Hcal$ be a Hilbert space and  $1\leq p\le n_o <\infty$.
 The function $z\mapsto \det_{n_o}(1-z A)$ is entire for every fixed
 $A\in \Scal_p(\Hcal)$
 and has the  representation
\begin{equation}
    \label{det-m-cn-eq}%
\det\nolimits_{n_o}(1-z A)=
 1+\sum_{n={n_o}}^\infty \frac{c_n(A)}{n!}\,(-1)^n\,  z^n,
\end{equation}
 where the coefficients $c_n(A)$ are defined by
  $$
 c_n(A):=
 \det\left(\begin{smallmatrix}
 b_1&n-1&0&\dots &0&0\\
 b_2&b_1&n-2&\dots& 0&0\\
 b_3&\ddots&\ddots &\ddots&\vdots&\vdots\\
 \vdots&\vdots&\vdots &\ddots&\ddots&\vdots\\
 b_{n-1}&b_{n-2}&b_{n-3}&\dots&b_1&1\\
 b_n&b_{n-1}&b_{n-2}&\dots&b_2&b_1
 \end{smallmatrix}\right)
 $$
   and
 $$
 b_n:=\left\{
 \begin{array}{ll}
 \trace A^n,&\mbox{ if } n\geq n_o,\\
 0,&\mbox{ otherwise.}
 \end{array}
 \right.$$
For $|z|$ sufficiently small one has
\begin{equation}
    \label{det-m-exptrace-eq}%
\det\nolimits_{n_o}(1-z A)\,=\, \exp\!\Big(-
\sum_{n=n_o}^\infty \frac{z^n}{n}\, \trace A^n\Big).
\end{equation}
 Let $(\lambda_j)_j$ be the eigenvalues of
$A\in \Scal_p(\Hcal)$, then one  has the Euler product
\begin{equation}
    \label{det-m-exptrace-euler-eq}%
 \det\nolimits_{n_o} (1-zA)\,=\,
 \prod_j\Big((1- z\lambda_j)\,\exp\!\big(\sum_{k=1}^{n_o-1}\frac{\lambda_j^k}{k}\, z^k\big)\Big)
 \,=\,
 \prod_j f_{n_o}(z\lambda_j),
\end{equation}
 where $\ds
 f_{n_o}(z)=(1- z)\,\exp\!\Big(\sum_{k=1}^{n_o-1}\frac{z^k}{k}\Big)
\stackrel{(\star)}{=}\exp\!\Big(-\sum_{k=n_o}^\infty \frac{z^k}{k}\Big)
 $.
 \qed
 \end{Lemma}
 The identity $(\star)$ can be obtained as  a consequence of the power series expansion of $\log(1-z)$.
The Euler product expansion (\ref{det-m-exptrace-euler-eq}) shows that the
zeros of $z\mapsto  \det\nolimits_{n_o} (1-zA)$ are in bijection with the eigenvalues of $A$.
\begin{Lemma}
    \label{gogokr-lemma-m2}%
Let $\Hcal$ be a Hilbert space.
Then  for any $n_o\in \N$ there exists  a constant $\Gamma_{n_o}>0$
such that for all $A\in \Scal_{n_o}(\Hcal)$ the estimates below hold.
One can choose  $\Gamma_1=1$.
\begin{equation}
\label{est-det-m}%
|\det\nolimits_{n_o}(1+A)|\,
\leq\,
\exp\!\big({\Gamma_{n_o}\|A^{n_o}\|_{\Scal_1(\Hcal)}}\big)
\, \leq\,
\exp\!\big({\Gamma_{n_o}\|A\|_{\Scal_{n_o}(\Hcal)}^{n_o}}\big).
\end{equation}
\end{Lemma}

\begin{proof}
The inequality
$ |\det\nolimits_{n_o}(1+A)|
\leq
 \exp\!\big(\Gamma_{n_o}\|A\|_{\Scal_{n_o}(\Hcal)}^{n_o}\big)$
can be found in   \cite[Thm.~XI.2.2]{GoGoKr00}. It based on
the estimate $|f_{n_o}(z)|\leq \exp\!\big(\Gamma_{n_o}\,|z|^{n_o}\big)$
for some constant $\Gamma_{n_o}>0$. Using (\ref{det-m-exptrace-euler-eq}) one derives from this the
slightly sharper estimate
\begin{eqnarray*}
 |\det\nolimits_{n_o}(1+A)|
&\leq& \exp\!\big(\Gamma_{n_o} \sum_j |\lambda_j(A)|^{n_o}\big)\\
&=& \exp\!\big(\Gamma_{n_o} \sum_j |\lambda_j(A^{n_o})|\big)\\
&\leq& \exp\!\big(\Gamma_{n_o} \| A^{n_o} \|_{\Scal_1(\Hcal)} \big).
\end{eqnarray*}
\end{proof}

The following criterion for the convergence of infinite products of regularized determinants
will turn out to be useful.

\begin{Lemma}
     \label{main-lemma}%
Let $(\Hcal_\nu)_{\nu\in \N}$   be a family of
Hilbert spaces. Fix $n_o\in \N$ and pick
$G_{\nu}\in \Scal_{n_o}(\Hcal_\nu)$
satisfying $
\sum_{\nu=0}^{\infty}\| G_{\nu}^{n_o}\|_{\Scal_{1}(\Hcal_\nu)}<\infty.
$
Then
$$
\prod_{\nu=0}^\infty\det\nolimits_{n_o} (1-z G_{\nu}).
$$
converges absolutely and locally uniformly to an entire function of $z$.
\end{Lemma}

\begin{proof}
Note that the
function $f_{n_o}$ from Lemma~\ref{gogokr-lemma-m}  is of the form
$f_{n_o}(z)-1=z^{n_o}\, h_{n_o}(z)$  for some entire function $h_{n_o}$.
Further,
 $c:=\sup_{\nu\in \N_0 }\|G_{\nu}\|$ is finite.
Hence   $|h_{n_o}(\lambda_j(z G_\nu))|
 \leq \sup_{|w|\leq |z|\, c} |h_{n_o}(w)|=:c_z
 $ for all eigenvalues $\lambda_j(zG_\nu)$ of $zG_\nu$.
 Now the hypothesis
 implies that
\begin{eqnarray*}
\sum_{\nu=0}^\infty\sum_j |f_{n_o}(\lambda_j(z G_{\nu}))-1|
&\leq &\sum_{\nu=0}^\infty\sum_j |\lambda_j(z G_{\nu})|^{n_o}\,
       |h_{n_o}(\lambda_j(z G_{\nu}))|\\
&\leq & c_z\, |z|^{n_o}\,\sum_{\nu=0}^\infty\sum_j |\lambda_j( G_{\nu}^{n_o})|\\
&\leq & c_z\, |z|^{n_o}\,\sum_{\nu=0}^\infty\sum_j|s_j( G_{\nu}^{n_o})|\\
&=& c_z\, |z|^{n_o}\,
   \sum_{\nu=0}^\infty  \|G_{\nu}^{n_o}\|_{\Scal_{1}(\Hcal_\nu)}
\end{eqnarray*}
is finite. Thus the infinite product
$\prod_{\nu=0}^\infty\prod_j f_{n_o}(\lambda_j(z G_\nu))$ converges, and by
Lemma~\ref{gogokr-lemma-m}
it is equal to
$\prod_{\nu=0}^\infty\det\nolimits_{n_o} (1-z G_\nu)$. This proves the claim.
\end{proof}

\begin{Prop}
 \label{det-ev-prop}
Let $\Hcal_1$ and $\Hcal_2$ be two Hilbert spaces and $n_o\in \N$. Pick
$A\in \Scal_{n_o}(\Hcal_1)$  and $B\in \Scal_{n_o}(\Hcal_2)$. If we denote the
eigenvalues of $B$ by  $\lambda_j(B)$, we have
$$
\det\nolimits_{n_o}\big(1-z \, A\otimes B\big)
=\prod_j \det\nolimits_{n_o}\big(1-z \lambda_j(B)\,A\big).
$$
\end{Prop}
\begin{proof}
For $|z|<\|A\|_{\Scal_{n_o}(\Hcal_1)}^{-1}\, \|B\|_{\Scal_{n_o}(\Hcal_2)}^{-1}$
the Lidskii Trace Theorem (\cite[Thm. IV.6.1]{GoGoKr00})
applied to the trace class operators $B^n$ ($n\geq {n_o}$) yields
\begin{eqnarray*}
\det\nolimits_{n_o}(1-z A\otimes B)
&\stackrel{(\ref{det-m-exptrace-eq})}{=}&\exp\!\Big(-\sum_{n={n_o}}^\infty
     \frac{z^n}{n}\, \trace (A\otimes B)^n\Big)\\
&=&\exp\!\Big(-\sum_{n={n_o}}^\infty \frac{z^n}{n}\, \trace A^n\ \trace B^n\Big)\\
&=&\exp\!\Big(-\sum_{n={n_o}}^\infty \frac{z^n}{n}\, \trace A^n\
   \sum_j\lambda_j(B)^n\Big)\\
&=&\exp\!\Big(- \sum_j\sum_{n={n_o}}^\infty \frac{z^n}{n}\,
    \lambda_j(B)^n\, \trace A^n\Big)\\
&=&\prod_j \exp\!\Big(-\sum_{n={n_o}}^\infty \frac{z^n}{n}\,
    \trace(\lambda_j(B)\, A)^n\Big)\\
&\stackrel{(\ref{det-m-exptrace-eq})}{=}&
\prod_j\det\nolimits_{n_o}\big(1-z \lambda_j(B)\,A\big).
\end{eqnarray*}
Here we used that for $z$ in the chosen range the inner double series converges
absolutely and locally uniformly.
By Lemma \ref{gogokr-lemma-m} the left hand side is an entire function in $z$.
Therefore analytic continuation shows that the identity  holds for all $z\in \C$ if we can show
that also the right hand side is an entire function in $z$.
But that follows from
Lemma~\ref{main-lemma} applied to the family $G_j:=  \lambda_j(B)\,A$.
\end{proof}

Let $A:\Hcal\rightarrow \Hcal$ be a trace class operator on a Hilbert space $\Hcal$
and
$ \wedge^rA:\bigwedge^r\Hcal\rightarrow\bigwedge^r\Hcal$
its $r$-fold exterior product. Then (cf. \cite{S77}) we have
\begin{equation}\label{multilinear-lemma}
\det(1-A)=\sum_{r=0}^{\dim \Hcal} (-1)^r \,\trace\! \wedge^r A
\end{equation}
and the
estimate
 \begin{equation}\label{multilinear-lemma-est}
 \| \wedge^r A \|_{\Scal_{1}( \wedge^r\Hcal)}
\leq {\textstyle\frac{1}{r!}}\,  \| A \|_{\Scal_{1}( \Hcal)}^r.
\end{equation}
For the special case of  a finite rank operator $B$ with spectrum
 $\lambda_1, \ldots, \lambda_d$
Proposition~\ref{det-ev-prop} implies  that
\begin{equation}\label{det-ev-cor}
\det\nolimits_{n_o}(1-z A\otimes \wedge^\nu B)
=
\prod_{\alpha\in \{0,1\}^d;\, |\alpha|=\nu}
 \det\nolimits_{n_o}(1-z \lambda^\alpha A),
\end{equation}
where for $\alpha\in\{0,1\}^d$ we set
 ${\lambda}^\alpha:=\prod_{j=1}^d\lambda_j^{ \alpha_j}$.
Approximating Schatten class operators by finite rank operators
one derives the following proposition:

\begin{Prop}\label{det-ev-cor-2}
For ${n_o}\in \N$ consider $A\in \Scal_{n_o}(\Hcal_1)$ and $B\in \Scal_{n_o}(\Hcal_2)$
with eigenvalues $\lambda_j:=\lambda_j(B)$. Then
\begin{enumerate}
\item[{\rm(i)}]
$$
\det\nolimits_{n_o}(1-z A\otimes \wedge^\nu B)
=\lim_{d\to\infty} \prod_{\alpha\in \{0,1\}^d;\, |\alpha|=\nu}
\det\nolimits_{n_o}(1-z \lambda^\alpha A).
$$
\item[{\rm(ii)}]  \label{main-prop}
$$
 \prod_{\nu=0}^\infty
 \det\nolimits_{n_o}(1-z  A\otimes \wedge^\nu B)
 =
 \lim_{d\to \infty} \prod_{\alpha\in \{0,1\}^d }
\det\nolimits_{n_o}(1-z \lambda^\alpha A).
 $$
\end{enumerate}
\end{Prop}

\begin{proof}
\begin{enumerate}
\item[(i)]
For any $d\in \N$ let $\pr_d\in \End(\Hcal_2)$ be the orthogonal projection onto the space
spanned by the first $d$ generalized eigenvectors of $B$ and
 $B_d:=\pr_d\circ B \circ\pr_d$.  Then   $(B_d)_{d\in \N}$ converges to $B $ in
 $\Scal_{n_o}(\Hcal_2)$ by \cite[Thm. IV 5.5]{GoGoKr00}.
Thus for any $C\in \Scal_{n_o}(\Hcal_1)$ the sequence
 $(C\otimes B_d)_{d\in \N}$ converges   to $C\otimes B$ in
 $\Scal_{n_o}(\Hcal_1\hat{\otimes}\Hcal_2)$.
The continuity of the regularized determinant
$X\mapsto\det\nolimits_{n_o}(1-X)$
now shows
$\lim_{d\to \infty }
\det\nolimits_{n_o}(1-z A\otimes B_d)=
\det\nolimits_{n_o}(1-z A\otimes  B)
$
for all $A\in \Scal_{n_o}(\Hcal_1)$, so that
(\ref{det-ev-cor}),
applied to the $B_d$,  implies the claim.
\item[(ii)]
  We know that
 \begin{eqnarray*}
\prod_{\alpha\in \{0,1\}^d }
\det\nolimits_{n_o}(1-z \lambda^\alpha A)
&=&
\prod_{\nu=0}^d \prod_{\alpha\in \{0,1\}^d;\, |\alpha|=\nu}
 \det\nolimits_{n_o}(1-z \lambda^\alpha A)
\\&=&
\prod_{\nu=0}^d
 \det\nolimits_{n_o}(1-z  A\otimes \wedge^\nu B_d)
\\&=&
\prod_{\nu=0}^\infty
 \det\nolimits_{n_o}(1-z  A\otimes \wedge^\nu B_d),
\end{eqnarray*}
since $\wedge^\nu C=0$ for any  $C\in \Mat(d,d;\C)$ and $\nu>d$.
In view of (i) it only remains to show that
$$
\lim_{d\to \infty}
\prod_{\nu=0}^\infty
 \det\nolimits_{n_o}(1-z  A\otimes \wedge^\nu B_d)
 =
\prod_{\nu=0}^\infty
\lim_{d\to \infty}
 \det\nolimits_{n_o}(1-z  A\otimes \wedge^\nu B_d).
$$
Expanding both sides in terms of the eigenvalues of $B$ we see that
it suffices to show
$$
\prod_{\nu=0}^\infty  \prod_{\alpha \in \{0,1\}^d;\, |\alpha|=\nu}
\!\! \det\nolimits_{n_o}(1-z \,\lambda^\alpha A)
\underset{d\to\infty}{\longrightarrow}
\prod_{\nu=0}^\infty   \prod_{\alpha \in \{0,1\}^\N;\, |\alpha|=\nu}
 \!\! \det\nolimits_{n_o}(1-z \,\lambda^\alpha A).
$$
For this kind of rearrangement it is enough to
 verify the summability condition  (\ref{main-thm-1-condition})
 from Lemma~\ref{main-lemma}, i.\!~e., the finiteness of
 $$
\sum_{\nu=0}^\infty \sum_{\alpha\in \{0,1\}^\N;\,|\alpha|=\nu}
\|z \,\lambda^\alpha\, A\|_{\Scal_{n_o}(\Hcal)}=
\|z \, A\|_{\Scal_{n_o}(\Hcal)}
\,
\sum_{\nu=0}^\infty \sum_{\alpha\in \{0,1\}^\N;\,|\alpha|=\nu}
|\lambda^\alpha|
.
 $$
 To show this, we recall that the eigenvalues of the trace class operator
$\wedge^\nu B$ are the $\lambda^\alpha$ with $|\alpha|=\nu$ and note that
\begin{eqnarray*}
\sum_{\nu=0}^\infty \sum_{\alpha\in\{0,1\}^\N;\,|\alpha|=\nu} |\lambda^\alpha|
&=& \sum_{\nu=0}^\infty  \sum_k |\lambda_k(\wedge^\nu B)|
\\
&\leq&\sum_{\nu=0}^\infty  \sum_k s_k(\wedge^\nu B)
\\
&=&\sum_{\nu=0}^\infty \|\wedge^\nu B\|_{\Scal_1(\wedge^\nu \Hcal_2)}
\\
&\stackrel{ (\ref{multilinear-lemma-est})}{ \leq}&\sum_{\nu=0}^\infty\frac{1}{\nu!}\,  \| B\|_{\Scal_1(\Hcal_2)}^\nu
<\infty.
\end{eqnarray*}
\end{enumerate}
\end{proof}

\begin{Lemma}{}
    \label{lemma-T-HS}%
Let $(F,\nu)$ be a measure space, $g:F\times F\rightarrow\C$ a
measurable function, and  $(S_x)_{x\in F}$  a measurable   family of
operators on a separable Hilbert space $\Hcal$.
 The formula
\begin{equation}
\label{t-schlange}
(T (f_1\otimes f_2))(\sigma):=\int_F g(x,\sigma)\,f_1(x)\, S_x f_2\, d\nu(x)
\end{equation}
defines a Hilbert-Schmidt operator $T:  L^2(F,d\nu)\hat{\otimes} \Hcal
\to  L^2(F,d\nu)\hat{\otimes} \Hcal$ iff
\begin{equation}
\label{t-schlange2}
\int_F \int_F|g(x,\sigma)|^2\, \,  d\nu(\sigma)\,
\|S_x \|^2_{\Scal_2( \Hcal)}\, d\nu(x)<\infty.
\end{equation}
 In this case $ T$ satisfies
$$
 \|T\|_{\Scal_2(L^2(F,d\nu)\hat{\otimes} \Hcal)}^2=
\int_F \int_F|g(x,\sigma)|^2\, \,  d\nu(\sigma)\,
\|S_x \|^2_{\Scal_2( \Hcal)}\, d\nu(x)
$$
and
$$
\trace  T^n=
\int_{F^n} \Big(\prod_{j=1}^{n-1} g(x_j,x_{j+1})\Big) \, g(x_n,x_1)\,\trace(S_{x_n}\circ \ldots\circ S_{x_1}) \, d\nu^n(x_1,\ldots,x_n)
$$
for  all $n\geq 2$. Moreover, for these $n$ we have
\begin{eqnarray*}
\| T^n\|^2_{\Scal_2( L^2(F,d\nu)\hat{\otimes} \Hcal)}
&=&
\int_F \int_{F^n}\big |\Big(\prod_{j=1}^{n-1} g(x_j,x_{j+1})\Big)\, g(x_n,\sigma)\big|^2\,\times\\
&&\quad
\times\, \|  S_{x_n}\circ \ldots\circ S_{x_1} \|^2_{\Scal_2(\Hcal)} \, d\nu^n(x_1,\ldots ,x_n)
  \, d\nu(\sigma)  .
  \end{eqnarray*}
\end{Lemma}

\begin{proof}
Suppose first that (\ref{t-schlange}) defines a Hilbert-Schmidt operator.
Fix orthonormal bases $(e_i)_{i\in \N}$, $(f_j)_{j\in \N}$  for $ L^2(F,d\nu)$ and $\Hcal$,
respectively. Then, by Parseval's identity, one has
\begin{eqnarray*}
\lefteqn{
\!\!\!\!
\|T\|^2_{\Scal_2( L^2(F,d\nu)\hat{\otimes} \Hcal)}
\,=\,
\sum_{i,j=1}^\infty \| T (e_i\otimes f_j)\|^2 }
\\
&=&
\sum_{i,j,k,l=1}^\infty\Big |\hermsp{ T(e_i\otimes f_j)}{e_k\otimes f_l}\Big|^2
\\
&=&
\sum_{i,j,k,l=1}^\infty \Big|
\int_F \int_F g(x,\sigma)\, e_i(x)\, \hermsp{S_x f_j}{f_l}\, d\nu(x)\, \overline{e_k(\sigma)}\, d\nu(\sigma)
\Big|^2
\\
&=&
\sum_{j,l=1}^\infty
\int_F \int_F \Big|g(x,\sigma)\,  \hermsp{S_x f_j}{f_l}\Big|^2\, d\nu(x)\,  d\nu(\sigma)
 \\
&=&
\int_F \int_F|g(x,\sigma)|^2\,
\sum_{j,l=1}^\infty   \Big|\hermsp{S_x f_j}{f_l}\Big|^2\, d\nu(x)\,  d\nu(\sigma)
\\
&=&
\int_F \int_F|g(x,\sigma)|^2\, \,  d\nu(\sigma)\,
\|S_x \|^2_{\Scal_2( \Hcal)}\, d\nu(x).
\end{eqnarray*}
Conversely, if (\ref{t-schlange2}) holds, we reverse this calculation and conclude
that not only the integral (\ref{t-schlange}) converges for almost all $\sigma$, but also that
it defines a Hilbert-Schmidt operator on $L^2(F,d\nu)\hat{\otimes} \Hcal$.

Now assume that $T$ is Hilbert-Schmidt. Then for $n\geq2$ the
operator $T^n$ is trace class and a simple induction argument shows
that
 \begin{eqnarray*}
 \lefteqn{
( T^n (e\otimes f))(\sigma)=
}\\
&=& \int_{F^n} \Big(\prod_{j=1}^{n-1} g(x_j,x_{j+1})\Big)\,
g(x_n,\sigma)\, e(x_1)\,S_{x_n}\circ \ldots\circ S_{x_1} f\,
d\nu^n(x_1,\ldots ,x_n) .
\end{eqnarray*}
By the first part of the proof the $S_{x_j}$ are Hilbert-Schmidt (for almost all $x_j$), hence the
compositions $S_{x_n}\circ \ldots\circ S_{x_1}$ are trace class.
Now the trace of $T^n$ can be calculated as follows
\begin{eqnarray*}
\trace T^{n}
&=&
\sum_{i,j=1}^\infty  \hermsp{ T^{n}(e_i\otimes f_j)}{e_i\otimes f_j}\\
&=&
\sum_{i,j=1}^\infty\int_F \int_{F^{n} }g(x_1,x_2) \cdots g(x_{n-1},x_n)\, g(x_n,\sigma)\,
\times\\
&& \phantom{\sum}
 \times\,\hermsp{S_{x_n}\circ \ldots \circ S_{x_1} f_j}{f_j}\, e_i(x_1)\,
 d\nu^n(x_1, \ldots  ,x_n)
\, \overline{e_i(\sigma)}\, d\nu(\sigma)\\
&=&
\sum_{i=1}^\infty \int_F \int_{F^{n} } g(x_1,x_2) \cdots g(x_{n-1},x_n)\, g(x_n,\sigma)\,
\times\\
&& \phantom{\sum}
 \times\, \trace (S_{x_n}\circ \ldots \circ S_{x_1}) \, e_i(x_1)\,
 d\nu^n(x_1, \ldots  ,x_n)
\, \overline{e_i(\sigma)}\, d\nu(\sigma).
\end{eqnarray*}
We claim that $\trace T^n$ can be rewritten as $\sum_{i=1}^\infty
\hermsp{ \Gcal_n e_i}{e_i} \,=\,\trace \Gcal_n $  with
\begin{eqnarray*}
(\Gcal_n f)(\sigma) & := &
  \int_{F^{n} } \Big(\prod_{j=1}^{n-1}g(x_j,x_{j+1})\Big) \, g(x_n, \sigma)\,
\times\\
&& \phantom{\sum}
 \times\, \trace (S_{x_n}\circ \ldots \circ S_{x_1}) \, f(x_1)\,
 d\nu^n(x_1, \ldots  ,x_n).
 \nonumber
\end{eqnarray*}
Note that (by Fourier expansion  and induction)
$$
\trace (S_n\circ \ldots \circ S_1)\ =
 \sum_{i_1,\ldots, i_n=1}^\infty \Big(\prod_{j=1}^{n-1} \hermsp{ S_j h_{i_j}}{h_{i_{j+1}}} \Big)\, \hermsp{ S_n h_{i_n}}{h_{i_{1}}}
$$
for any orthonormal   basis $(h_i)_{i\in \N}$ for $\Hcal$ and
Hilbert-Schmidt operators $S_i$ on $\Hcal$. Setting
$$
(\Gcal_{i,j} f)(\sigma)  :=
  \int_{F } g(x, \sigma)\,\hermsp{ S_x h_i}{h_{j}}   \, f(x)\, d\nu(x)
  $$
  for $i,\, j\in \N$,
we can
rewrite  $\Gcal_n$ as
 \begin{eqnarray*}
(\Gcal_n f)(\sigma) & = &
  \int_{F^{n} } \Big(\prod_{j=1}^{n-1}g(x_j,x_{j+1})\Big) \, g(x_n, \sigma)\,
\times\\
&& \phantom{\sum}
 \times\,
\sum_{i_1,\ldots, i_n=1}^\infty \prod_{j=1}^n \hermsp{ S_{x_i} h_{i_j}}{h_{i_{j+1}}}
 \, f(x_1)\,
 d\nu^n(x_1, \ldots  ,x_n)
   \\& = &
\sum_{i_1,\ldots, i_n=1}^\infty
(\Gcal_{i_n,i_1}\circ\Gcal_{i_{n-1}, i_n}\circ \ldots \circ \Gcal_{i_1,i_2} f)(\sigma).
  \end{eqnarray*}
The identity
\begin{eqnarray}\ \ \
\sum_{i, i=1}^\infty
    \| \Gcal_{i, j } \|_{\Scal_2(L^2(F,d\nu))}^2\!\!\!
 &=&\!\!\!
\sum_{i, j=1}^\infty \int_{F^2} | g(x,y)|^2\,|\hermsp{ S_x h_i}{h_{j}}|^2\, d\nu(x)\, d\nu(y)
\label{identity-gij}
 \\&=& \!\!\!\int_{F^2} | g(x,y)|^2\,
\sum_{i, j=1}^\infty  |\hermsp{ S_x h_i}{h_{j}}|^2\, d\nu(x)\, d\nu(y)
\nonumber
 \\&=&\!\!\! \int_{F^2} | g(x,y)|^2\,
   \| S_x\|^2_{\Scal_2(\Hcal)} \, d\nu(x)\, d\nu(y)
   \nonumber
      \end{eqnarray}
implies  that the $\Gcal_{i,j}$
are Hilbert-Schmidt operators on $L^2(F,d\nu)$. Therefore,
for each $(i_1,\ldots, i_n)\in \N^n$ the  integral operator
$ \Gcal_{i_n,i_1}\circ\Gcal_{i_{n-1}, i_n}\circ \ldots \circ \Gcal_{i_1,i_2}$
is trace class and  by \cite[Ex. X. 1.18]{Ka66}  its trace can be obtained by
integrating the integral kernel along the diagonal.
If  $\Gcal_n$ is trace class, we have
\begin{eqnarray*}
\lefteqn{
\trace \Gcal_n
 \,=\,
\sum_{i_1,\ldots, i_n=1}^\infty \trace
(\Gcal_{i_n,i_1}\circ\Gcal_{i_{n-1}, i_n}\circ \ldots \circ \Gcal_{i_1,i_2} )}
\\&=&
\int_{F^n} \Big(\prod_{j=1}^{n-1}g(x_j,x_{j+1})\Big) \, g(x_n,x_1)\,  \trace{(S_{x_n}\circ \ldots\circ S_{x_1} )}\,  d\nu^n(x_1,\ldots, x_n)
\\&=&
\trace T^{n}.
\end{eqnarray*}
Thus, to prove the claim it suffices to show that
$
\sum_{i_1,\ldots, i_n=1}^\infty
\Gcal_{i_n,i_1}\circ 
 \ldots \circ \Gcal_{i_1,i_2}$ converges in $\Scal_1\big(L^2(F, d\nu)\big)$.
 Using the technical Lemma~\ref{cycl-lemma} below, we obtain the estimate
\begin{eqnarray*}
\| \Gcal_n\|_{\Scal_1\big(L^2(F, d\nu)\big)} &\leq& \sum_{i_1,
\ldots, i_n=1}^\infty \prod_{j=1}^n \|\Gcal_{i_j,
i_{j+1}}\|_{\Scal_2\big(L^2(F, d\nu)\big)}
\\&
\leq&
\Big(\sum_{i, j=1}^\infty
\|\Gcal_{i_j, i_{j+1}}\|_{\Scal_2\big(L^2(F, d\nu)\big)}^2\Big)^{n/2}
 \\&
\stackrel{(\ref{identity-gij})}{=}&
 \Big(\int_F \int_F  |g(x, y)|^2 \,\| S_x\|^2_{\Scal_2(\Hcal)}
\, d\nu(x)\, d\nu(y) \Big)^{n/2},
\end{eqnarray*}
which proves the claim.
 To conclude the proof of the lemma one verifies the formula for
$\|T^n\|^2_{\Scal_2( L^2(F,d\nu)\hat{\otimes} \Hcal)}$
for $n\ge 2$, which can be done similarly as in the case $n=1$.
\end{proof}

If $\nu$ is a finite measure on $F$, $g:F^2\to\C$ is bounded,   and $\int_F \|S_x\|_{\Scal_2( \Hcal)}^2\, d\nu(x)$ is finite,
Lemma~\ref{lemma-T-HS}  shows   that
the associated  operator $T$  is Hilbert-Schmidt.

 \begin{Lemma}
    \label{cycl-lemma}
    Let $n\geq 2$  and suppose that the functions
  $a_k:\N\times\N\to \C$ satisfy
$
\sum_{i,j=1}^\infty |a_k(i,j)|^2 <\infty
$
 for $k=1, \ldots, n$.   Then
$$
\Big|\sum_{i_1, \ldots, i_n=1}^\infty
\prod_{k=1}^n a_k(i_k, i_{k+1}) \Big|
\leq
\prod_{k=1}^n
\Big(
\sum_{i,j=1}^\infty |a_k(i,j)|^2
\Big)^{1/2}
$$
using the convention that $i_{n+1}=i_1$.
\end{Lemma}
\begin{proof}
We proceed by induction. The case  $n=2$ follows from the estimate
\begin{eqnarray*}
\Big|\sum_{i_1, i_2=1}^\infty
a_1(i_1,i_2)\, a_2(i_2, i_1)  \Big|
\!&\leq &\!
\sum_{i_1=1}^\infty \Big(\sum_{ i_2=1}^\infty
|a_1(i_1,i_2)|^2\Big)^{1/2}\, \Big(\sum_{i_2=1}^\infty |a_2(i_2, i_1)|^2  \Big)^{1/2}
\\&\leq   &\!
\prod_{k=1}^2
\Big(
\sum_{i,j=1}^\infty |a_k(i,j)|^2
\Big)^{1/2}.
\end{eqnarray*}
To do the induction step  consider
$\tilde{a}_n(i,j):= \sum_{m=1}^\infty |a_n(i, m)\, a_{n+1}(m,j)|$. Then
\begin{eqnarray*}
 \sum_{i, j=1}^\infty |\tilde{a}_n(i,j)|^2
 & =&
 \sum_{i, j=1}^\infty\Big(
\sum_{m=1}^\infty |a_n(i, m)\, a_{n+1}(m,j)| \Big)^2
\\&\leq &
 \Big(\sum_{i, m=1}^\infty |a_n(i,m)|^2\Big) \,
 \Big(\sum_{ m,j=1}^\infty    |a_{n+1}(m,j)|^2\Big)  ,
\end{eqnarray*}
and induction yields
\begin{eqnarray*}
\Big|\sum_{i_1, \ldots,i_{n+1}=1}^\infty
\prod_{k=1}^{n+1} a_k(i_k, i_{k+1}) \Big|
&\leq &
\sum_{i_1, \ldots, i_{n}=1}^\infty
\Big|\prod_{k=1}^{n-1} a_k(i_k, i_{k+1}) \Big| \, \tilde{a}_n(i_n, i_1)\\
&\leq&
\prod_{k=1}^{n-1}
\Big(
\sum_{i,j=1}^\infty |a_k(i,j)|^2
\Big)^{1/2}\,
\Big( \sum_{i, j=1}^\infty |\tilde{a}_n(i,j)|^2
\Big)^{1/2}
\\
&\leq &
\prod_{k=1}^{n+1}
\Big(
\sum_{i,j=1}^\infty |a_k(i,j)|^2
\Big)^{1/2}.
\end{eqnarray*}
\end{proof}

\section{Meromorphic Continuation}
\label{MeromCont}

The following theorem
is an analog of a
result of D. Mayer
(see \cite[Thm. 7.17]{M91}) proven there in the context of
generalized Perron-Frobenius operators associated
with expanding maps.

\begin{Theorem}
    \label{main-thm-1}
Let $(\Hcal_\nu)_{\nu\in \N_0}$   be a
family of Hilbert spaces. Fix $n_o\in \N$ and pick $G_{\nu}\in
\Scal_{n_o}(\Hcal_\nu)$
 such that
\begin{equation}\label{main-thm-1-condition}
 \sum_{\nu=0}^{\infty}\| G_{\nu}^{n_o}\|_{\Scal_{1}(\Hcal_\nu)}<\infty.
\end{equation}
Let $(\zn)_{n\in \N}$ be a sequence in $\C$ such that $ \zn
=\sum_{\nu=0}^{\infty} (-1)^\nu \,\trace  G_{\nu}^n $ for all $n\geq
n_o$. Then the dynamical zeta function $\zeta_R$ associated with
$(\zn)_{n\in \N}$ admits a meromorphic continuation to the entire
plane. It is given by the formula
$$
\zeta_R(z)=
\exp\!\Big(\sum_{n=1}^{n_o-1} \frac{z^n}{n}\,\zn\Big)
\,
\prod_{\nu=0}^{\infty }\Big(\det\nolimits_{n_o} (1-z G_{\nu})\Big)^{(-1)^{\nu+1}}.
$$
\end{Theorem}
\begin{proof}
We treat the case of finitely many non-zero operators  first, say $G_l=0$ for all $l\geq k$.
For
$|z|<\min\{\|(G_{\nu})^{n_o}\|_{\Scal_1(\Hcal_\nu)}^{-1}\,|\,\nu=0,\ldots,k\}$, using
Lemma~\ref{gogokr-lemma-m}, one calculates
 \begin{eqnarray*}
\zeta_R(z)
&=&\exp\!\Big(\sum_{n=1}^{n_o-1} \frac{z^n}{n}\,\zn\Big)
    \,
    \exp\!\Big(\sum_{n=n_o}^\infty \sum_{\nu=0}^{k}(-1)^\nu\,
    \frac{z^n}{n}\,\trace G_{\nu}^n\Big)\\
&=&\exp\!\Big(\sum_{n=1}^{n_o-1} \frac{z^n}{n}\,\zn\Big)
   \,
   \prod_{\nu=0}^{k}\exp\!\Big(\sum_{n=n_o}^\infty  \frac{z^n}{n}\,\trace G_{\nu}^n\Big)^{(-1)^\nu} \\
&=&\exp\!\Big(\sum_{n=1}^{n_o-1} \frac{z^n}{n}\,\zn\Big)
   \,
   \prod_{\nu=0}^{k}\Big(\det\nolimits_{n_o} (1-z G_{\nu})\Big)^{(-1)^{\nu+1}}
\end{eqnarray*}
and obtains that this is a finite product of meromorphic functions.

We turn the general case. By (\ref{main-thm-1-condition}) the sequence
$\| G_{\nu}^{n_o}\|_{\Scal_{1}(\Hcal_\nu)}$ tends to zero as $\nu\to\infty$,
so the minimum
$\min\{\|G_{\nu}^{n_o}\|_{\Scal_1(\Hcal_\nu)}^{-1}\,|\,\nu\in \No\}>0$ exists.
Using the convergence criterion from Lemma~\ref{main-lemma}
we see that the quotient of infinite products
$$
\prod_{\nu=0}^\infty\Big(\det\nolimits_{n_o} (1-z G_{\nu})\Big)^{(-1)^{\nu+1}}=
\frac{
\prod_{\nu=0}^\infty\det\nolimits_{n_o} (1-z G_{2\nu+1})
}{
\prod_{\nu=0}^\infty\det\nolimits_{n_o} (1-z G_{2\nu})
}
$$
converges absolutely and locally uniformly for all $z\in\C$.
\end{proof}

\begin{Corollary}
    \label{main-thm-2}
Given two Hilbert spaces $\Hcal$ and $\Hcal_o$, fix $n_o\in \N$,  and consider $G\in \Scal_{n_o}(\Hcal)$
and $\Lambda\in \Scal_{n_o}(\Hcal_o)$.
Let $(\zn)_{n\in \N}$ be a sequence in $\C$ such that
$
\zn= \det(1-\Lambda^n)\,\trace G^n
$
for $n\geq n_o$.
Then the dynamical
zeta function $\zeta_R$ associated with $(\zn)_{n\in \N}$
admits a meromorphic continuation to the entire plane. It
is given by the formula
$$
\zeta_R(z)
=
\exp\!\Big(\sum_{n=1}^{n_o-1} \frac{z^n}{n}\,\zn\Big)
\,
\prod_{\nu=0}^{\dim \Hcal_o}\Big(\det\nolimits_{n_o}
\big(1-z\, G\otimes \wedge^\nu\Lambda\big)\Big)^{(-1)^{\nu+1}}.
$$
\end{Corollary}
\begin{proof}
Since $\wedge^r(A^n)=(\wedge^r A)^n$ and $\trace A\ \trace B=\trace(A\otimes B)$
for all trace class operators $A$ and $B$, the identity (\ref{multilinear-lemma})
implies  that
\begin{eqnarray*}
Z_n
&=&\det(1-\Lambda^n)\, \trace G^n\\
&=&\sum_{\nu=0}^{\dim \Hcal_o} (-1)^\nu \,\trace\! (\wedge^\nu \Lambda^n)\, \trace G^n\\
&=&\sum_{\nu=0}^{\dim \Hcal_o} (-1)^\nu \,\trace G_{\nu}^n
\end{eqnarray*}
with  $G_{\nu}:= G\otimes \wedge^\nu\Lambda $ on
$\Hcal^{(\nu)}:=\Hcal\hat{\otimes}  \wedge^\nu\Hcal_o $ for $n\ge n_o$.
Note here that the estimate (\ref{multilinear-lemma-est})
provides  the summability condition  (\ref{main-thm-1-condition}).
 In fact,
\begin{eqnarray*}
\sum_{\nu=0}^{\dim \Hcal_o} \| (G_{\nu})^{n_o}\|_{\Scal_{1}(\Hcal^{(\nu)})}
&=&
\sum_{\nu=0}^{\dim \Hcal_o} \| G^{n_o}\otimes \wedge^\nu\Lambda^{n_o} \|_{\Scal_1(\Hcal^{(\nu)})}
\\
&=&
\| G^{n_o}\|_{\Scal_1(\Hcal)}\,\sum_{\nu=0}^{\dim \Hcal_o}  \|
\wedge^\nu\Lambda^{n_o} \|_{\Scal_1( \wedge^\nu\Hcal_o)}
\\
&\leq&
\| G^{n_o}\|_{\Scal_1(\Hcal)}\,\sum_{\nu=0}^{\infty } \frac{1}{\nu!}\,  \| \Lambda^{n_o} \|_{\Scal_{1}( \Hcal_o)}^\nu
<\infty.
\end{eqnarray*}
Now we can apply  Theorem~\ref{main-thm-1} to finish the proof.
\end{proof}

Consider the case of
a dynamical zeta function in the presence of the
dynamical trace formula of the type $Z_n=\det(1-\Lambda^n)\, \trace G^n$.
We obtain a better understanding the zeros and poles of the zeta function
calculating the regularized
determinants as infinite products involving the eigenvalues of $\Lambda$.
In fact, combining Proposition \ref{det-ev-cor-2} with Corollary \ref{main-thm-2} yields
the following spectral interpretation of the poles and zeros of $\zeta_R$.

\begin{Theorem} \label{main-thm-3}
Given two Hilbert spaces $\Hcal$ and $\Hcal_o$ fix $n_o\in \N$ and consider $G\in \Scal_{n_o}(\Hcal)$
and $\Lambda\in \Scal_{n_o}(\Hcal_o)$.
Let $(\zn)_{n\in \N}$ be a sequence in $\C$ such that
$
\zn= \det(1-\Lambda^n)\,\trace G^n
$
for $n\geq n_o$.
Denote the eigenvalues of $\Lambda$ by  $(\lambda_i)_{i\in\N}$,
repeated according to multiplicity.
For $\alpha\in\{0,1\}^d$ set
 ${\lambda}^\alpha:=\prod_{\nu=1}^d\lambda_\nu^{ \alpha_\nu}$.
 Then the meromorphic continuation of
 $\zeta_R$ is given by the formula
$$
\zeta_R(z) = \exp\!\Big(\sum_{n=1}^{n_o-1} \frac{z^n}{n}\,\zn\Big)
\, \lim_{d\rightarrow \infty}\prod_{\alpha\in\{0,1\}^d} \Big(
\det\nolimits_{n_o}(1-z{\lambda}^\alpha\,  G)\Big)^{ (-1)^{|\alpha|+1}}.
$$
 \qed
 \end{Theorem}

\section{Composition Operators on Fock Spaces}
    \label{fock-sec}
We start by briefly recalling some basic properties of reproducing
kernel spaces of holomorphic functions on not necessarily finite
dimensional manifolds. Our basic reference for this material is
\cite{Ne00}, although we choose a different normalization.

Let $\Hcal\subset\C^E$ be a Hilbert space consisting of complex
valued functions on a set $E$. The space $\Hcal$ is called a
\textit{reproducing kernel Hilbert space (RKHS)}, if for each $x\in
E$ the evaluation functional $ \ev_x:\Hcal\rightarrow \C,\ f\mapsto
f(x) $ is continuous. A function $k:E\times E\rightarrow \C$ is
called a \textit{reproducing kernel} for $\Hcal$, if for all $y\in
E$ the function $k_y:=k(\,\cdot\,,y):E\rightarrow \C$ belongs to
$\Hcal$ and if for all $f\in \Hcal$, $y\in E$ we have $
f(y)=\hermsp{ f }{ k_y }_\Hcal $. Recall that a function $p:E\times
E\rightarrow \C$ is \textit{positive definite}, if $ \sum_{k,l=1}^n
\overline{a_k}\,a_l\,p(x_k,x_l)\geq 0 $ for all $n\in\N$, $a_j\in
\C$, $x_j\in E$ ($j=1,\ldots,n$).
These concepts are connected by the following fact,
 cf. \cite[I.1]{Ne00}:
If $\Hcal\subset\C^E$ is a  RKHS,
then the function
$k:E\times E\rightarrow \C$ defined by
$k(x,y):=\ev_x\circ \ev_y^*$ is a positive definite reproducing kernel for $\Hcal$.
Moreover, for all
$f\in \Hcal$, $ x\in E$ one has
\begin{equation}\label{standard-rkhs-cor}
|f(x)|\leq \|f\|\, \sqrt{k(x,x)}.
\end{equation}
\label{fock-ex}%
Since the span of the kernel functions $k_w$ ($w\in E$) is dense in $\Hcal$,
a bounded operator $T$ on   $\Hcal$    is uniquely determined by  its ``integral kernel''
\begin{equation}
\label{integral-kernel-t}
k_T(z,w):=
(T k_w)(z)=
\hermsp{T k_w}{k_z}  .
\end{equation}
The integral kernel of the adjoint $T^*$ of $T$ is obtained from
$$
k_{T^*}(z,w)=
\hermsp{T^* k_w}{k_z}  =
\hermsp{k_w}{T k_z}  =\overline{\hermsp{T k_z}{k_w}}   = \overline{k_T(w,z)}.
$$

Recall that the \textit{Bargmann-Fock space}
     $\Fcal(\C^m)$  is defined as the  space of entire functions
     $F:\C^m\rightarrow \C$ with
$$
\|f\|^2_{\Fcal(\C^m)}:=\int_{\C^m}|f(z)|^2\, \exp(-\pi \|z\|^2)\,dz<\infty
$$
where $dz$ denotes Lebesgue measure on $\C^m$.
It is a RKHS with  reproducing kernel $\ds k(z,w)= \exp\big(\pi\langle z|w\rangle\big).$

 Let  $(\Hcal,\langle\cdot|\cdot\rangle)$  be a separable Hilbert space, then  the map
$
k:\Hcal\times \Hcal\to\C,\,(z,w)\mapsto \exp(\pi \langle z|w\rangle )$
is a positive definite kernel, see  \cite[I.2.2]{Ne00}.
  One defines the \textit{(symmetric) Fock space}
 to be the unique reproducing kernel Hilbert space
 $\Fcal(\Hcal)\subset \C^{\Hcal}$ associated with this kernel.
Since the reproducing kernel
is  holomorphic in the first and  anti-holomorphic in the second variable,
 \cite[Prop. A.III.10]{Ne00} shows that
$\Fcal(\Hcal )\subset \Ocal(\Hcal )$. In particular,  $\Fcal(\Hcal )$ consists of continuous functions.
Let $(e_i)_{i\in I}$ be an orthonormal  basis for $\Hcal$.  Then the monomials
$
\zeta_\alpha(z)= \big(\frac{\pi^\alpha}{\alpha!}\big)^{1/2}\, \prod_i \langle z | e_i\rangle^{\alpha_i}$
($\alpha\in (\N_0)^I$)
form an orthonormal basis for the  Fock space $\Fcal(\Hcal)$,
where we use the standard multiindex notations:
$\alpha!:=\prod_i \alpha_i!$ and     $ |\alpha|:=\sum_i \alpha_i$.
We will use  composition operators  to view the $\Fcal(\C^m)$ as subspaces of $\Fcal(\Hcal)$. In this context we note the
following elementary lemma.

 \begin{Lemma}
    \label{lemma-1}%
    Let $\Hcal_1,\, \Hcal_2$ be separable Hilbert spaces and $A:\Hcal_1\to\Hcal_2$ a linear operator
    with $\|A\|\leq 1$.
    Define $C_A: \Fcal(\Hcal_2) \to \Fcal(\Hcal_1), \ f\mapsto f\circ A$. Then
    $C_A$ is continuous with
    $\|C_A\|\leq 1 $ and  $(C_A)^*=C_{A^*}$.
If $A$ is surjective, then  $C_A$ is injective.
    In particular, if $AA^*=\id$, then  $C_A$ is an isometric embedding.
 \end{Lemma}
 \begin{proof}
  We   use   the reproducing kernel property for the reproducing kernels
  $k^{(i)}$ of $\Fcal(\Hcal_i)$ ($i=1,2$).
   Let $x\in \Hcal_2, \, y\in \Hcal_1$. Then
$$
\langle C_A k_x^{(2)} | k_y^{(1)} \rangle =
(k_x^{(2)}\circ A)(y)= k^{(2)}(Ay,x)=e^{\pi \langle Ay|x\rangle_{2}}
=e^{\pi \langle y|A^* x\rangle_1}
$$
and
$$
e^{\pi \langle y|A^* x\rangle_1}
= k^{(1)}(y,A^* x)
= (C_{A^*} k_x^{(1)})(y)= \langle C_{A^*} k_x^{(1)} | k_y^{(2)}\rangle.
$$
 This  shows that $C_{A}$ and $C_{A^*}$ are adjoint.
   \begin{eqnarray*}
  \|C_A k_x^{(2)}\|^2
  &=&
\langle C_A k_x^{(2)} |C_A k_x^{(2)}\rangle
=
\langle C_{A^*} C_A k_x^{(2)} | k_x^{(2)}\rangle
=
\langle C_{AA^*} k_x^{(2) }| k_x^{(2)}\rangle
\\
&=&
k^{(2)}(AA^* x,x)
=
k^{(2)}(A^* x,A^*x)
=\|k^{(2)}_{A^* x}\|^2,
\end{eqnarray*}
so
$$
  \|C_A k_x^{(2)}\|
=\|k^{(2)}_{A^* x}\|
= \exp({\textstyle \frac{1}{2} }  \|A^* x\|^2)
\leq
\exp({\textstyle \frac{1}{2} }   \|x\|^2)
=\|k^{(2)}_{ x}\|
$$
implies the claim.
\end{proof}

Using  Lemma~\ref{lemma-1} we obtain
\begin{enumerate}
\item
If $P:\Hcal\to \Hcal$  is  a
projection, then $C_P\in \End(\Fcal(\Hcal))$ is a projection.
\item
If $P:\Hcal\to \Hcal$  is  self-adjoint,  then $C_P\in \End(\Fcal(\Hcal))$ is self-adjoint.
\item
As a consequence of (i) and (ii) we see that if  $P:\Hcal\to \Hcal$  is an  orthogonal projection,  then $C_P\in \End(\Fcal(\Hcal))$ is an  orthogonal projection.
\item
    Let $P:\Hcal\to \Hcal$  be an orthogonal projection. For  $p:\Hcal\to P\Hcal, \, z\mapsto Pz$ we have
      $P=p^*p$ and      $pp^*=\id\nolimits_{P\Hcal}$. Hence
  $C_{p}:\Fcal(P\Hcal)\rightarrow\Fcal(\Hcal)$
   is an isometric embedding.
   \item With the identification $\Fcal(P\Hcal)\cong C_{p}(\Fcal(P\Hcal))\subset  \Fcal(\Hcal)$
    we can view $\Fcal(P\Hcal)$ as a subspace of $ \Fcal(\Hcal)$.
Moreover, $\Fcal(P\Hcal)$ has a reproducing kernel, namely the kernel of $\Fcal(\Hcal)$ restricted to
 $P\Hcal\times P\Hcal$.
\item  $C_{p*}:\Fcal(\Hcal)\rightarrow \Fcal(P\Hcal)$
is the adjoint of $C_{p}:\Fcal(P\Hcal)\rightarrow\Fcal(\Hcal) $,
 $C_P=C_p(C_p)^*$, and
$1=\|C_{p}\|=\|C_{p*}\|$.
\end{enumerate}

 \begin{Lemma}
    \label{lemma-prn*}%
    Let $P_n:\Hcal\to \Hcal$ be a sequence of orthogonal projections
 converging  to the identity in the strong operator topology.
 Then the   sequence $C_{P_n} \in \End(\Fcal(\Hcal))$
 of orthogonal projections
 converges  in the strong operator topology
 to the identity  on $\Fcal(\Hcal)$ as $n\rightarrow \infty$.
\end{Lemma}\begin{proof}
For all  $z\in \Hcal$ one has $P_n z \longrightarrow  z$ as $n\to \infty$.
Since $ \Fcal(\Hcal)\subset  \Ccal(\Hcal)$, we have for all $f\in  \Fcal(\Hcal)$
$$
C_{P_n} f(z)
= f(P_n z ) \longrightarrow f(z).
$$
Using the reproducing kernel property, this can be rewritten as
$$
\hermsp{C_{P_n} f}{k_z} \longrightarrow \hermsp{f}{k_z}
$$
for all $z\in \Hcal$. Since the functions $k_z$ ($z\in \Hcal$) form a total subset of $ \Fcal(\Hcal)$, this  implies weak operator convergence which   on a
  Hilbert space  coincides with the strong operator convergence.
\end{proof}

Parts of the following proposition can be found in \cite[II.]{Ma88}.

\begin{Prop}    \label{fock-thm*}%
    Let $P_n:\Hcal\to \Hcal$ be an ascending sequence
 of orthogonal projections with $n$-dimensional range
 converging  in the strong operator topology to the identity, i.\!~e., $P_n\Hcal\subset P_{n+1}\Hcal$. Set
 $\pr_n:\Hcal\to \Hcal_n:=P_n\Hcal, \, z\mapsto P_nz$.
 A function $f$ belongs to the Fock space  $\Fcal(\Hcal)$, defined as the RKHS with reproducing kernel
$k(z,w)=\exp(\pi \langle z|w\rangle)$,
if and only if the following three conditions hold:
 \begin{enumerate}
 \item
 $f:\Hcal\rightarrow \C$ is continuous,
 \item
for all $m\in \N$ the map $f\circ \pr_m^*:\Hcal_m\rightarrow \C$ is analytic, and
 \item
 $\ds \sup_{m\in \N} \int_{\Hcal_m} |f\circ \pr_m^*(z)|^2\, \exp(-\pi \|z\|^2)\,
 dz<\infty$.
\end{enumerate}
 In this case, this supremum equals
 $$ \|f\|_{\Fcal(\Hcal)}^2=\lim_{m\in \N} \int_{\Hcal_m} |f\circ \pr_m^*(z)|^2\, \exp(-\pi \|z\|^2)\,  dz
 .
 $$
\end{Prop}

\begin{proof}
We already observed that
$\Fcal(\Hcal )\subset \Ocal(\Hcal )$, hence (i) and (ii) hold for $f\in \Fcal(\Hcal)$.
For $m\in \N$ we introduce
\begin{eqnarray*}
c_m(f)&:=& \int_{\Hcal_m} |f\circ \pr_m^*(z)|^2\, \exp(-\pi \|z\|^2)\,  dz
\\&=&\|C_{\pr_m^*}f\|^2_{\Fcal(\Hcal_m)} =
\|C_{\pr_m}C_{\pr_m^*}f\|^2_{\Fcal(\Hcal )}=
\|C_{P_m}f\|^2_{\Fcal(\Hcal )}.
\end{eqnarray*}
  By Lemma~\ref{lemma-prn*} one gets
$$ \lim_{m\to \infty} c_m(f)= \lim_{m\to \infty} \|C_{P_m}f\|^2_{\Fcal(\Hcal )}
=\| \lim_{m\to \infty}  C_{P_m}f\|^2_{\Fcal(\Hcal )}=\|f\|^2_{\Fcal(\Hcal )}
.
$$
Suppose, conversely that  $f:\Hcal \rightarrow \C$ satisfies (i) - (iii). The bounded sequence
$ \big(C_{P_m}f\big)_{m\in \N}$ in $ \Ccal(\Hcal )$ has a weakly convergent subsequence with limit $f$.
Hence  $f\in \Fcal(\Hcal )$  and
 $$
 \|f\|_{\Fcal(\Hcal )}^2 = \lim_{m\to \infty } c_m(f) \leq \sup_{m\in \N} c_m(f).
 $$
 Viewing
 $\Fcal(\Hcal_m)$ as a subspace of  $\Fcal(\Hcal_{m+1})$,
 Parseval's identity implies that the sequence
$c_m(f)
$
 indexed by $m\in \N$ is monotonically increasing.
\end{proof}


Let $E$ be a set and $V$ a space of complex valued functions on $E$.
A \textit{(weighted} or \textit{generalized) composition operator}
is an operator $T:V\to V$ of the form
$$
(Tf)(z)=\phi(z)\, (f\circ \psi)(z),
$$
where
$\phi:E\rightarrow \C$, $\psi:E\rightarrow E$ are given functions.
If the multiplication part is trivial, i.\!~e., $\phi\equiv 1$, then $T$ is called a
\textit{(classical) composition operator}.

Let  $E,\,F$ be non-empty sets.  Let $\phi_x:E\rightarrow \C$, $\psi_x:E\rightarrow E$ for each $x\in F$, and
 $T_x:\C^E\rightarrow \C^E$, $ (T_xf)(z):= \phi_x(z)\, (f\circ\psi_x)(z).$
Then a simple induction argument yields
the composition law
\begin{equation}\label{comp-op-iterates}%
(T_{x_n}\circ \ldots \circ T_{x_1}f)(z)=
\prod_{k=1}^n (\phi_{x_k}\circ \psi_{x_{k+1}}\circ\ldots\circ\psi_{x_n})(z)
\,
(f\circ\psi_{x_1}\circ\ldots\circ \psi_{x_n})(z).
\end{equation}

\begin{Remark}
    \label{kontrak-1-prop}%
Let $0<q<1$ and $\psi:X\rightarrow X$ be a function on a normed space $(X,\|\,\cdot\,\|)$ with
$
\|\psi(z)-\psi(w)\|\leq q\,\|z-w\|$
for all $z,\, w\in  X$. Then $\psi$ is called a   contraction.
Set
\begin{equation}
    \label{standard-r-estimate}%
r_\psi:=\frac{\|\psi(0)\|}{1-q}.
\end{equation}
Let $  r> r_\psi$ and  $|z|\leq r.$ Then the estimate
$$
\|\psi(z)\|\leq \|\psi(z)-\psi(0)\|+\|\psi(0)\|
\leq q\,\|z\|+\|\psi(0)\|<r
$$
shows that the set  $K_r:=\big\{z\in X\ \big|\  \|z\|\leq r\big\}$ satisfies
$\psi(K_r)\subset K_{qr+\|\psi(0)\|}\subset K_r$.
Let $\psi^{(m)}:=\psi\circ \ldots\circ \psi$ ($m$-times)
be the $m$-th iterate of $\psi$.
Then for any $z\in X$ and
$m\geq n_0\geq \frac{\ln(\frac{r-r_\psi}{\|z\|})}{\ln q}$ we have
$$
\|\psi^{(m)}(z)\|
\leq
\|\psi^{(m)}(z)-\psi^{(m)}(0)\|
+\|\psi^{(m)}(0)\|
\leq q^m \,\|z\|+r_\psi
\leq r,
$$
since $0\in K_{r_\psi}$  implies  $\psi^{(m)}(0)\in K_{r_\psi}$.
\qed
\end{Remark}

\begin{Remark}
    \label{kontrak-ganz-prop}%
    Let $(X,\|\,\cdot\,\|)$ be a normed space,
    $\psi: X\to X$ be a contraction in the sense of Proposition~\ref{kontrak-1-prop}, and
$\phi:X\rightarrow \C$  a continuous function. Let
$r>r_\psi$ with $r_\psi$ as in~(\ref{standard-r-estimate}), and
 $T$ be the weighted  composition operator
$$
T:\Ccal(K_r)\to \Ccal(K_r),\
(Tf)(z)= \phi(z)\, (f\circ\psi)(z).
$$
\begin{enumerate}
\item
Let $g\in \Ccal(K_r)$, then
 $T g$ belongs to $\Ccal(K_{\delta r})$ for some  $\delta>1$.
In fact, if  $z\in K_{\delta r}$ and $\delta< \frac{r-\|\psi(0)\|}{rq}$,
then $|\psi(z)|\leq q\,r\delta+\|\psi(0)\|< r$.
Since  $\frac{r-\|\psi(0)\|}{rq}>1$, we may
choose $\delta >1$.

\item Every  eigenfunction of $T$  for a non-zero eigenvalue belongs to $ \Ccal(X)$.
To see this, let
$f\in \Ccal(K_r)$ be an eigenfunction of $T$  for a non-zero eigenvalue $\rho$. Hence by
iterating relation~(i) $n$-times we get  $f=\rho^{-n}T^n f \in  \Ccal(K_{\delta^n r})$ for
some  $\delta>1$. Since  $X=\bigcup_{r>0}K_r$,
    we conclude that  $f\in \Ccal(X)$.
\qed
\end{enumerate}
 \end{Remark}

\section{A Trace Formula for Composition Operators}
    \label{traces-sec}

Let $U\subset \C^k$ be an open bounded complex domain.
Let $A^\infty(U)$ denote the  space of holomorphic functions on $U$ which are continuous on  the closure
$\overline{U} $ of $U$. Clearly,  $A^\infty(U)$  is a Banach space with respect to
the supremum norm.

 The following theorem, due to D. Ruelle (\cite{Ru76}, see also \cite[Appendix B]{M80a},
\cite{M80b} for the infinite dimensional case)
 is based on a fixed point formula of Atiyah and Bott (cf. \cite{AB67}).
\begin{Theorem}
            \label{atiyah-bott}%
Let $U\subset \C^k$ be an open bounded complex domain.  Let $\phi:U\to \C$ and $ \psi:U\to U$ be  holomorphic functions with continuous extensions to
$\overline{U} $  and, moreover,
  $\psi(\overline{U})\subset U$.
Then $\psi$ has a unique fixed point  $z^*\in U$ and the weighted
composition operator
$$
T :A^\infty(U)\rightarrow A^\infty(U),  \ (T f)(z)= \phi(z) \, (f \circ \psi)(z)
$$
is nuclear of order zero
 with trace given by  the Atiyah-Bott type fixed point formula
 $$
\trace_{\! A^\infty(U)\,} T= \frac{\phi(z^*)}{\det(1-\psi^\prime (z^*))}.
$$
\qed
\end{Theorem}

\begin{Lemma}
    \label{kontrak2-lemma}%
Let  $\psi:\C^m\to \C^m$  and  $\phi:\C^m\to\C$ be entire  functions, and  $\psi$ a contraction in the sense of Proposition~\ref{kontrak-1-prop}.
Let
$r>r_\psi$ with $r_\psi$ as in~(\ref{standard-r-estimate}) and
$T:A^\infty(B(0;r))\rightarrow A^\infty(B(0;r)) $
be the composition operator  acting via
$$(Tf)(z)=\phi(z)\, (f\circ\psi)(z).
$$
Let
 $f$  an
 eigenfunction of  $T$
for  a non-zero eigenvalue $\rho$. Then $f$ is  entire and
 there exist $c_1,\, c_2 >0$ such that for all $z\in \C^m$
$$
|f(z)|\leq \|z\|^{-c_1\ln\rho}\, \sup_{|w|\leq r}|f(w)| \,\max_{t\in [0,2\pi]}
 | \phi(e^{it}z))|^{c_2 \ln\|z\| }
.$$
Moreover, if  $A^2(U):=\Ocal(U)\cap L^2(U,dz)$ denotes the  Bergman space, then
$$
 \trace_{\!A^\infty(U)} T= \trace_{\!A^2(U)} T
$$
for all  $\psi$-invariant bounded domains $U\subset \C^m$.
\end{Lemma}
\begin{proof}
Let $f$ be an eigenfunction of $T $
for a non-zero eigenvalue $\rho$. For $n\in \N$ we have
$f=\rho^{-n}\,T^n f$ which
 by (\ref{comp-op-iterates}) is given as
$$
f(z)=\rho^{-n} \prod_{k=0}^{n-1}(\phi\circ \psi^{(k)} )(z)
\,(f\circ\psi^{(n)})(z),
$$
where $\psi^{(k)}$ is the $k$-th iterate of $\psi$.
As in Remark~\ref{kontrak-ganz-prop} (ii)
one shows  that   $f$ is entire,
thus belongs to $A^2(U)$ for all bounded domains $U\subset \C^m$.
Hence every eigenvalue of $T|_{A^\infty(U)} $ belongs to the spectrum of $T|_{A^2(U)} $, thus by Lidskii's Trace Theorem the traces coincide.

Let  $z\in \C^d$ with $\|z\|> r_\psi$ and
 $n(z):=\Big\lceil \frac{\ln(\frac{r-r_\psi}{\|z\|})}{\ln q}\Big\rceil$.
 One can find  constants $c_1,\,c_2>0$ such that
 $c_1\ln\|z\|\leq n(z)\leq c_2\ln\|z\|$ for all $\|z\|> r_\psi$.
 Remark~\ref{kontrak-1-prop} implies that
  $\|\psi^{(n(z))}(z)\|\leq r$, and hence
\begin{eqnarray*}
|f(z)| &=&
|\rho|^{-n(z)} \ \Big|\prod_{k=0}^{n(z)-1}(\phi\circ \psi^{(k)} )(z)\Big|
\,\big|(f\circ\psi^{(n(z))}(z)\big|
\\
&\leq& |\rho|^{-n(z)}\, \sup_{|w|\leq r}|f(w)| \,\sup_{\|w\|\leq \|z\|}
\Big|\prod_{k=0}^{n(z)-1}(\phi\circ \psi^{(k)} )(w)\Big|
\\
&\leq& |\rho|^{-n(z)}\, \sup_{|w|\leq r}|f(w)| \,\sup_{\|w\|\leq \|z\|}|\phi(w)|^{n(z)}
\\
&\leq& |\rho|^{-c_1\ln \|z\|}\, \sup_{|w|\leq r}|f(w)| \,
\sup_{\|w\|\leq \|z\|} | \phi(w)|^{c_2 \ln\|z\| }.
\end{eqnarray*}
By the maximum principle  we know that the supremum
$\sup_{\|w\|\leq \|z\|}|\phi(w)|$ is attained for some $w$ with $\|w\|=\|z\|$.
\end{proof}

\begin{Theorem}
    \label{spurformel-fock-thm}%
Let
$b\in \C^m$,   $\A\in \Gl(m;\C)$ with
$\|\A\|<1$, and $\phi:\C^m\to \C$ an entire function which  can be estimated
by
$|\phi(z)|\leq c \exp(a\,\|z\|)$
for some constants $a,\,c>0$. Let
$T$ be the composition operator given by
$$
(Tf)(z)=\phi(z)\, f(\A z+ b)
.$$
Then $T:\Fcal(\C^m)\rightarrow \Fcal(\C^m)$
is a trace class
operator with
$$\trace_{\!\Fcal(\C^m)} T
=
\trace_{\!A^\infty(B(0;r))} T=
\frac{\phi((1-\A)^{-1}b)}{\det(1-\A)}
$$
for all $B(0;r):=\big\{z\in \C^m\ \big|\  \|z\|<r\}$ with $r>\frac{\|b\|}{1-\|\A\|}$.
\end{Theorem}

\begin{proof}
The affine map  $\psi(z)=\A z+b$ is a contraction with  $q=\|\A\|<1$
and  $r_\psi=\frac{\|b\|}{1-\|\A\|}$.
We claim that the Fock space $\Fcal(\C^m)$  is a $T$-invariant Hilbert subspace of
$A^\infty(B(0;r))$ for any  $r>r_\psi=\frac{\|b\|}{1-\|\A\|}$.
In fact, for $f\in \Fcal(\C^m)$ the  standard estimate (\ref{standard-rkhs-cor}) yields
\begin{eqnarray*}
\|T f\|^2
&=&
\int_{\C^m} \Big|
\phi(z)\,f(\A z+b)\Big|^2 e^{-\pi\|z\|^2}\,dz
\\
&\leq&
c^2
\int_{\C^m} e^{2a\|z\|}\,|f(\A z+b)|^2 \,e^{-\pi\|z\|^2}\,dz
\\
&\leq &c^2\, \|f\|^2\,
\int_{\C^m}e^{2a\|z\|}\, e^{\pi\|\A z+b\|^2} \,e^{-\pi\|z\|^2}\,dz
\\
&\leq& \|f\|^2\,
\Big( C+
\int_{\C^m\setminus B(0;r)}  e^{2a\|z\|}\,  e^{-\pi(1-\|\A\|^2)\|z\|^2}\,dz
\Big)
<\infty.
\end{eqnarray*}
Thus $T|_{\Fcal(\C^m)}$ is a nuclear map  on a Hilbert space, and hence of trace class.
Let $f\in A^\infty(B(0;r))$ be an eigenfunction of $T$ corresponding  to a non-zero eigenvalue $\rho$.
By Lemma~\ref{kontrak2-lemma}  the eigenfunction $f$  satisfies the estimate
$$
|f(z)|^2\,\exp(-\pi\|z\|^2)\leq
 \|z\|^{-c_1\ln\rho}\,  \,\exp\big((a\, \|z\|\,+ \ln c)c_2\, \ln\|z\|\big)
\,\exp(-\pi\|z\|^2).
$$
This upper bound  is Lebesgue-integrable on $\C^m$,
and thus $f$ belongs to
$\Fcal(\C^m)$. This shows that  every non-zero eigenvalue of $T|_{A^\infty(B(0;r))}$
is an eigenvalue of
$T|_{\Fcal(\C^m)}$, hence the traces coincide and by Theorem~\ref{atiyah-bott} they have the stated value.
\end{proof}
Let $T$ be a trace class operator on a Hilbert space $\Hcal\subset L^2(Z,dm)$ with reproducing kernel $k$. Then by general theory  the trace is given by integrating the integral kernel (\ref{integral-kernel-t}) along the diagonal,
$$\trace T= \int_Z  k_T(z,z)\, dm(z)=\int_Z  (T k_z)(z)\, dm(z).
$$
 Thus Theorem~\ref{spurformel-fock-thm} yields the non-trivial integral identity
 \begin{equation}
\frac{\phi((1-\A)^{-1}b)}{\det(1-\A)}
=\int_{\C^m} \phi(z)\, e^{\pi\langle \A z+b|z\rangle}\, e^{-\pi\|z\|^2}\, dz.
\end{equation}

Let $\Hcal$ be a separable Hilbert space
and  $\Fcal(\Hcal) $   the associated  Fock space. Fix $a,\, b\in \Hcal$, and $\A\in \End(\Hcal)$.
Consider the  (possibly unbounded) composition operator
 \begin{equation}
    \label{kab-def}
\Kcal_{a,b,\A}:\Fcal(\Hcal)\to \Fcal(\Hcal),\
(\Kcal_{a,b,\A} f)(z)=e^{\pi \langle z|a\rangle}\,  f(\A z+b).
 \end{equation}
 If  $\Hcal$ is  finite dimensional and
 $\|\A\|<1$, then by
 Theorem~\ref{spurformel-fock-thm}
the operator $ \Kcal_{ a,b,\A}$ is trace class, hence compact.

 \begin{Prop}
    \label{kab-positive-prop}%
Let $(\Hcal, \hermsp{\cdot}{\cdot})$ be a finite dimensional Hilbert space.
\begin{enumerate}
\item
Let  $a,\,b\in \Hcal,\ \A\in \End(\Hcal)$ with $\|\A\|<1$, and
 $
  \Kcal_{ a,b,\A}\in \End( \Fcal(\Hcal) )
  $  be the corresponding composition operator
(\ref{kab-def}).
Then $(\Kcal_{a,b,\A})^\star =\Kcal_{b,a,\A^\star}$ and
 $\Kcal_{ a,b,\A}$ is selfadjoint if and only if $\A$ is selfadjoint and $a=b$.
 \item If
$\A$  is positive, then
 $\Kcal_{ b,b,\A}$ is positive and trace class with
$$
\trace \Kcal_{\beta,\beta,\A}
=\|\Kcal_{\beta,\beta,\A}\|_{\Scal_1(\Fcal(\Hcal))}=
 \frac{
\exp\!\big( \pi \|(1-\A)^{-1/2}\beta\|^2 \big)}{
\det(1-\A)}.
$$
 \item
 Let
 $a_i,\,b_i\in \Hcal$, $\A_i\in \End(\Hcal)$ with  $\|\A_i\|<1$ ($i=1,\,2$), then
$$
\Kcal_{a_1,b_1,\A_1}\Kcal_{a_2,b_2,\A_2}
=
e^{\pi \langle b_1|a_2\rangle}\,
\Kcal_{a_1+\A_1^\star a_2,\A_2 b_1+b_2,\A_2\A_1}.
$$
\end{enumerate}
 \end{Prop}

 \begin{proof}
The operator  $  \Kcal_{ a,b,\A}$  has  integral kernel
$$
k_ {\Kcal_{ a,b,\A}}(z,w)=  e^{\pi \langle z|a\rangle}\, e^{\pi\langle \A z+b| w\rangle}=
\exp\!\big(\pi ( \langle z|a\rangle+\langle b| w\rangle+ \langle \A z| w\rangle)\big),
$$
from which one easily gets the integral kernel
$$
k_ {(\Kcal_{ a,b,\A})^*}(z,w)=\overline{k_{\Kcal_{a,b,\A}}(w,z)}=
\exp\!\big(\pi \langle z| b\rangle+\pi\langle a+ \A^\star  z |w\rangle\big)
$$
of its adjoint $(\Kcal_{a,b,\A})^\star $. Similarly  one confirms that
$$
k_{\Kcal_{b,a,\A^\star}}(z,w)=
\exp\!\big(\pi \langle z| b\rangle+\pi\langle a+ \A^\star  z |w\rangle\big),
$$
i.\!~e., $(\Kcal_{a,b,\A})^\star =\Kcal_{b,a,\A^\star}$.
Hence $\Kcal$ is selfadjoint
iff
$\A$ is selfadjoint and $a=b$.
A compact selfadjoint operator is positive iff all its eigenvalues are positive.
By an argument given in \cite[III]{M80b}
we know the spectrum of  $\Kcal_{ b,b,\A}$: Set
${\phi}(z):= e^{\pi \langle z|b\rangle}$, $\psi(z):=\A z+b$, $z^*=(1-\A)^{-1}b$. Then
the spectrum  $\spec( \Kcal_{ b,b,\A}) $ of $\Kcal_{ b,b,\A}$ is contained in the set
\begin{eqnarray*}
\sigma^*&=&\{0\}\cup
  \big\{
  \phi(z^*)\,
\mu_{i_1}\cdot\ldots\cdot\mu_{i_k}\,|\, k\in \No,\,\mu_{j}\in \spec(\psi^\prime(z^*))\big\}
\\
&=&\{0\}\cup
\big\{e^{\pi \langle (1-\A)^{-1}b|b\rangle}
\mu_{i_1}\cdot\ldots\cdot\mu_{i_k}\,|\, k\in \No,\,\mu_{j}\in \spec(\A)\big\}.
\end{eqnarray*}
If $\A$ is positive,  then
all eigenvalues of $\Kcal_{ b,b,\A}$ are necessarily positive.
Hence the operator $\Kcal_{ b,b,\A}$
is positive and thus the trace norm $\|\Kcal_{\beta,\beta,\A}\|_{\Scal_1(\Fcal(\Hcal))}$
of $\Kcal_{\beta,\beta,\A}$ is
equal to its  trace, which is given by Theorem~\ref{spurformel-fock-thm}.
 A routine calculation confirms the composition law.
\end{proof}

\begin{Lemma}
    \label{trace-norm-and-square-root-lemma}%
Let $(\Hcal, \hermsp{\cdot}{\cdot})$ be a finite dimensional Hilbert space,  $\A\in \End(\Hcal)$
with $\|\A\|<1$,  and $a,\, b\in \Hcal$. Set
$$\Lambda=\sqrt{\A\A^\star},\ \beta=(1+\sqrt{\A\A^\star})^{-1}(\A a+b),
\,\gamma=\exp\!\Big(\frac{\pi}{2}(\|a\|^2-\|\beta\|^2)\Big).
$$
 Let
  $\Kcal:=\Kcal_{a,b,\A}$ and $ K:=\gamma\,\Kcal_{\beta,\beta,\Lambda}\in \End( \Fcal(\Hcal) )
  $ be the  corresponding composition operators (\ref{kab-def}).
  Then  $K=|\Kcal|=\sqrt{\Kcal^\star\Kcal}$, and
\begin{eqnarray*}
\|\Kcal\|_{\Scal_1(\Fcal(\Hcal))}
&=&
\frac{\gamma\,\exp\!\big( \pi \|(1-\Lambda)^{-1/2}\beta\|^2 \big) }{\det(1-\Lambda)}\\
&=&
\frac{\exp\!\big(
   \frac{\pi}{2}\|a\|^2+\frac{\pi}{2}\|(1-{\A\A^\star})^{-1/2}(\A a+b)\|^2\big)}{\det(1-|\A|)}.
\end{eqnarray*}
 \end{Lemma}

\begin{proof}
The composition operator
$K\in \End( \Fcal(\Hcal) )
 $ is  positive by Prop.~\ref{kab-positive-prop} (ii).
For all $f\in \Fcal(\Hcal)$ one has by Prop.~\ref{kab-positive-prop} (iii)
 \begin{eqnarray*}
(K^2 f)(z)&=&
\gamma^2\,e^{\pi\|\beta\|^2} \,
e^{\pi \langle z|(1+\Lambda^\star)\beta\rangle}\, f(\Lambda^2 z+(1+\Lambda)\beta),\\
(\Kcal^\star\Kcal f)(z)
&=&
e^{\pi\|a\|^2 }
e^{\pi\langle  z|\B a+b\rangle}\, f(\A\A^\star z+\A a+b).
\end{eqnarray*}
  For $\Lambda,\, \beta,\,\gamma$ chosen as above
 we get $K^2= \Kcal^\star \Kcal$, hence
$K=|\Kcal|=\sqrt{\Kcal^\star\Kcal}$
by the uniqueness of the operator square root.
 The trace norm of
$\Kcal$  is equal to  the trace of $K$, which is given  by  Prop.~\ref{kab-positive-prop} (ii).
In view of
\begin{eqnarray*}
\lefteqn{
2 \|(1-\Lambda )^{-1/2}\beta \|^2 -\|\beta \|^2\, =
}\\
&=& \!\!\!
2 \|(1-\sqrt{\A\A^\star})^{-1/2}(1+\sqrt{\A\A^\star})^{-1}(\A a+b)\|^2
-\|(1+\sqrt{\A\A^\star})^{-1}(\A a+b)\|^2
\\&=& \!\!\!
\langle
(2(1-\sqrt{\A\A^\star})^{-1}-1)
(1+\sqrt{\A\A^\star})^{-1}(\A a+b)|(1+\sqrt{\A\A^\star})^{-1}(\A a+b)\rangle
\\&=& \!\!\!
\langle (1+\sqrt{\A\A^\star})(1-\sqrt{\A\A^\star})^{-1}
(1+\sqrt{\A\A^\star})^{-1}(\A a+b)|(1+\sqrt{\A\A^\star})^{-1}(\A a+b)\rangle
\\&=& \!\!\!
\langle (1-\sqrt{\A\A^\star})^{-1}
(1+\sqrt{\A\A^\star})^{-1}(\A a+b)|\A a+b\rangle
\\&=& \!\!\!
\langle (1-{\A\A^\star})^{-1}(\A a+b)|\A a+b\rangle
\\&=& \!\!\!
\|(1-{\A\A^\star})^{-1/2}(\A a+b)\|^2
\end{eqnarray*}
we obtain
$$
\gamma \,
\exp\!\big( \pi \|(1-\Lambda )^{-1/2}\beta \|^2 \big)
=\exp\!\big(\frac{\pi}{2}\|a\|^2+\frac{\pi}{2}\|(1-{\A\A^\star})^{-1/2}(\A a+b)\|^2\big).
$$
 \end{proof}

Similarly to the proof of Lemma~\ref{trace-norm-and-square-root-lemma} one shows that
$ \sqrt{\Kcal\Kcal^\star}$ is given by
 $K':= \gamma'\,\Kcal_{\beta',\beta',\Lambda'}\in \End( \Fcal(\Hcal) )$ with
 $$\Lambda':=|\A|=\sqrt{\A^\star\A},\ \beta':=(1+|\A|)^{-1}(\A^\star b+a),
\,\gamma':=\exp\!\Big(\frac{\pi}{2}(\|b\|^2-\|\beta_1\|^2)\Big).
$$

\begin{Theorem}
    \label{spurformel-fockinfty-thm}%
 Let $a,\,b\in  \Hcal$,  $\A\in \Scal_1(\Hcal)$ with  $\|\A\|<1$, and
 consider the weighted composition operator
$$\Kcal: \Fcal(\Hcal) \to  \Fcal(\Hcal), \
(\Kcal f)(z)=e^{\pi \langle z|a\rangle}\,  f(\A z+b)
.$$
 \begin{enumerate}
\item
The operator  $K:=|\Kcal|:=\sqrt{\Kcal^\star\Kcal}:\Fcal(\Hcal) \to  \Fcal(\Hcal)$ is
given by
$(K f)(z):=\gamma\, e^{\pi \langle z|\beta\rangle}\, f(\Lambda z+\beta),
$
where
$$\Lambda=\sqrt{\A\A^\star},\ \beta=(1+\sqrt{\A\A^\star})^{-1}(\A a+b),
\,\gamma=\exp\!\Big(\frac{\pi}{2}(\|a\|^2-\|\beta\|^2)\Big)
$$
\item The operator $\Kcal$ is trace class with
$$
\|\Kcal\|_{\Scal_1(\Fcal(\Hcal))}=\trace K
=\frac{\exp\!\Big(\frac{\pi}{2}\|a\|^2+\frac{\pi}{2}\|(1-{\A\A^\star})^{-1/2}(\A a+b)\|^2\Big)}{\det(1-|\A|)}
 $$
and
$$
\trace \Kcal=\frac{\exp\!\big( \pi \langle (1-\A)^{-1}b|a\rangle \big)}{\det(1-\A)}.
$$
\end{enumerate}
\end{Theorem}

\begin{proof}
As in Lemma~\ref{trace-norm-and-square-root-lemma} one gets $K=|\Kcal|=\sqrt{\Kcal^\star\Kcal}$.
It remains to
show that the trace norm of $\Kcal$, i.\!~e., the trace  of $K$ is finite.
    Let $P_n:\Hcal\to \Hcal$ be an ascending sequence
 of orthogonal projections with $n$-dimensional range
 converging to the identity  in the strong operator topology. Set
 $\pr_n:\Hcal\to \Hcal_n:=P_n\Hcal, \, z\mapsto P_nz$.
 Fix $m\in \N$ and consider
 $K_m: = C_{\pr_m^*} K C_{\pr_m}\in \End( \Fcal(\Hcal_m) )
 $, which  by the composition law
acts via
$$
(K_m  f)(z)=\gamma\,
e^{\pi \langle z|\pr_m \beta\rangle}\, f(\pr_m\Lambda \pr_m^* z+\pr_m \beta)
=
e^{\pi \langle z|\beta_m\rangle}\, f(\Lambda_m z+\beta_m)
$$
with $\beta_m=\pr_m^*\beta$ and $\Lambda_m=\pr_m\Lambda \pr_m^*$.
 For any $v\in \Fcal(\Hcal_m)$ we have
$$
\hermsp{ K_m v }{ v}_{\Fcal(\Hcal_m)}
 =
 \hermsp{ C_{\pr_m^*} K C_{\pr_m} v  }{ v }_{\Fcal(\Hcal_m)}
 =
 \hermsp{ K v }{ v }_{\Fcal(\Hcal)}
 .$$
 Let $\{v_\alpha\,| \, \alpha \in I_m\}$
be an orthonormal basis for  $\Fcal(\Hcal_m)$, hence
$$
\sum_{\alpha\in I_m}\hermsp{ K v_\alpha }{ v_\alpha }
=
\sum_{\alpha\in I_m} \hermsp{ K_m\zeta_\alpha }{ \zeta_\alpha }_{\Fcal(\Hcal_m)}
=\trace K_m .
$$
The left hand side is monotonically increasing as $m\to \infty$. If it is bounded, it has limit
$\trace K$.
By the Atiyah-Bott formula from  Prop.~\ref{kab-positive-prop} (ii)      we have
$$
\trace K_m
=
\gamma\,\frac{\exp\!\Big( \pi \|(1-\Lambda_m)^{-1/2}\beta_m\|^2_{\Hcal_m} \Big)}{\det_{\Hcal_m}(1-\Lambda_m)}
=
\gamma\,\frac{\exp\!\Big( \pi \|(1-\Lambda_m)^{-1/2}\beta_m\|^2_{\Hcal} \Big)}{\det_{\Hcal}(1-\Lambda_m)}
$$
identifying $\Lambda_m \in \End( \Hcal_m)$ with $
\pr_m^*\Lambda_m\pr_m
  =P_m \Lambda P_m\in \End(\Hcal)$.
Lemma~\ref{lemma-prn*}, together with \cite[Thm. IV 5.5]{GoGoKr00}, shows
that the pointwise convergence  $P_m \Lambda P_m\rightarrow \Lambda$ is in trace norm.
Therefore the following limit exists:
\begin{eqnarray*}
\lim_{m\rightarrow\infty}
\trace K_m
&=&
\lim_{m\rightarrow\infty}\gamma\,
 \frac{\exp\!\Big( \pi \|(1-P_m \Lambda P_m)^{-1/2}\beta_m\|_\Hcal^2 \Big)}{\det_\Hcal(1-P_m\Lambda P_m)}
\\
&=&\gamma\,\frac{\exp\!\big( \pi \|(1-\Lambda)^{-1/2}\beta\|_\Hcal^2 \big)}{\det_\Hcal(1-\Lambda)}
<\infty.
\end{eqnarray*}
Thus $K$ and $\Kcal$  are trace class.
By Lemma~\ref{lemma-prn*}  and \cite[Thm. IV 5.5]{GoGoKr00} the sequence of  trace class operators
 $  C_{P_m} \Kcal C_{P_m}
 $ converges to $\Kcal$ in trace norm, hence the Atiyah-Bott fixed point formula
 from Proposition \ref{kab-positive-prop} allows to calculate
\begin{eqnarray*}
\trace \Kcal
&=&
\lim_{m\rightarrow\infty}
\trace C_{\pr_m}C_{\iota_m} \Kcal C_{\pr_m}C_{\iota_m}
=
\lim_{m\rightarrow\infty}
\trace    \Kcal_m
\\
&=&
\lim_{m\rightarrow\infty}\frac{\exp\!\big( \pi \hermsp{ (1-\A_m)^{-1}b_m }{ a_m }_{\Hcal_m} \big)}{\det_{\Hcal_m}(1-\A_m)}\,
\\
&=&\frac{\exp\!\big( \pi \hermsp{ (1-\A)^{-1}b }{ a }  \big)}{\det_\Hcal(1-\A)},
\end{eqnarray*}
where we view
$\A_m,{a}_m$, and ${b}_m$ as operator on, respectively vectors in, $\Hcal$.
A similar approximation argument, combined with Lemma~\ref{trace-norm-and-square-root-lemma}, also yields
$$
\trace K=
\frac{\exp\!\Big(\frac{\pi}{2}\|a\|^2+\frac{\pi}{2}\|(1-{\A\A^\star})^{-1/2}(\A a+b)\|^2\Big)}{\det_\Hcal(1-|\A|)}.
$$
\end{proof}

As a corollary we obtain an exact formula for the Hilbert-Schmidt
norm of a weighted composition operator of the form
$\Kcal_{a,b,\A}$ which could also be obtained directly as a consequence of
an identity on Gaussian integrals. We first consider the general situation: Let
$\Hcal\subset L^2(Z,dm)$ be a Hilbert space with reproducing kernel
$k$. Consider the composition operator
$$
(Tf)(z)=\phi(z)\, (f\circ \psi)(z),
$$
where
$\phi:Z\rightarrow \C$, $\psi:Z\rightarrow Z$ are fixed functions.
Then the Hilbert-Schmidt norm of $T$ is equal to
$$\int_Z |\phi(z)|^2\, k(\psi(z),\psi(z))\, dm(z).
$$
In fact: $T$ has the integral kernel $k_T(z,w)=\phi(z)\, k(\psi(z),w)$ and
\begin{eqnarray*}
\|T\|^2_{\Scal_2(\Hcal)} &=&
\int_Z \int_Z |k_T(z,w)|^2\, dm(w)\, dm(z)
\\&=&
\int_Z  |\phi(z)|^2\, \int_Z  k(w,\psi(z))\, \overline{k(w,\psi(z))}\,  dm(w)\, dm(z)
\\&=&
\int_Z |\phi(z)|^2\, k(\psi(z),\psi(z))\, dm(z)
\end{eqnarray*}
by the reproducing kernel property.

\begin{Corollary}
    \label{comp-positive-infty-cor}%
 Let $a,\,b\in  \Hcal$,  $\A\in \Scal_2(\Hcal)$ with  $\|\A\|<1$.
 Then the weighted composition operator
$$\Kcal: \Fcal(\Hcal) \to  \Fcal(\Hcal), \
(\Kcal f)(z)=e^{\pi \langle z|a\rangle}\,  f(\A z+b)
$$
is Hilbert-Schmidt with
$$\|\Kcal\|^2_{\Scal_2( \Fcal(\Hcal))}=
 \frac{\exp\!\Big( \pi\|a\|^2 + \pi \|   (1-\A \A^*)^{-1/2}(\A a+b)\|^2 \Big) }{\det(1-\A\A^*)}.
$$
\end{Corollary}

\begin{proof}
We use that
$
\|\Kcal\|^2_{\Scal_2( \Fcal(\Hcal))}=\trace \Kcal\Kcal^*
$
together with Theorem~\ref{spurformel-fockinfty-thm}.
By the arguments given in the proof of  Lemma~\ref{kab-positive-prop}
one has
$$
(\Kcal^\star\Kcal f)(z)=
e^{\pi\|a\|^2 }
e^{\pi\langle  z|\A a+b\rangle}\, f(\A\A^\star z+\A a+b).
$$
We note that $\|\A\A^\star\|=\|\A\|^2<1$ and $\A\A^\star\in  \Scal_1(\Hcal)$.
Hence the assertion follows from Theorem~\ref{spurformel-fockinfty-thm}.
\end{proof}

\section{Dynamical Trace Formulae for Spin Chains}
\label{dyntraceformulae}

Lemma~\ref{lemma-1} allows to describe the partition function of a non-interacting
matrix subshift as the trace of an operator, which is then called a {\it transfer operator}.
In the special case of a finite alphabet $F$ this result is well-known.

\begin{Prop}
\label{non-interacting-case-cor}
 Let  $(F,\nu,\A,0) $ be a matrix subshift with vanishing interaction.
 Then for $n\geq 2$ the integral  operator
$$
\Gcal_\A : L^2(F,d\nu)\rightarrow  L^2(F,d\nu), \
(\Gcal_\A f)(x)=
\int_F \A_{\sigma,x}\,f(\sigma)\, d\nu(\sigma)
$$
associated with the transition matrix $\A$
satisfies the dynamical trace formula
$$
\zn(0)= \nu^n\big(\rho_n(\{\xi\in \Omega_\A\mid \tau^n\xi=\xi\})\big)  =
\trace \Gcal_\A^n.
$$
\end{Prop}

\begin{proof}
The operator $ \Gcal_\A$ can be seen
as $T$ in Lemma~\ref{lemma-1}, where all the operators $S_x:=\id:\C\to\C$ are trivial.
Hence $ \Gcal_\A$ is Hilbert-Schmidt and
the traces of its iterates are given by Lemma~\ref{lemma-1}. Comparison of the formulae with
(\ref{partition-fct}) now proves the claim.
\end{proof}

 Let  $(F,\nu,\A,A) $ be a matrix subshift, where the interaction
 $A$ is of the following form:
Let $\B\in \Scal_p(\Hcal)$ for some $p<\infty$ with $\rhospec(\B)<1$.
For $x\in F$
assume that one has $a_x,\, b_x\in  \Hcal$,  $q_x\in \C$ such that
for $\xi=(x_j)_{j\in \N}\in \Omega_\A$
\begin{equation}
    \label{a-def}
A(\xi )= q_{x_1} + \pi \sum_{j=2}^\infty \hermsp{\B^{j-2}b_{x_{j}}}{a_{x_1}}.
\end{equation}
Here we assume absolut convergence of the series.
We interprete $A(\xi)$
as the sum of the  potential term $q_{x_1}$
and the sum of  two-body interactions between the first particle and the particles at positions
$i\in \N_{>1}$.

Let $\B\in \Scal_p(\Hcal)$ for some $p<\infty$ with $\rhospec(\B)<1$, then
 one can choose $n\in \N$
 large enough such that
 $\B^{n}$  has operator norm  less than one and is trace class, since the spectral radius can be  characterized via
$$
\rhospec(\B)=
\max\big\{|z|\in \R\,|\, z\in \spec(\B)\big\}=
\lim_{k\to\infty}\sqrt[k]{\|\B^k\|}
$$
and $\B^n \in \Scal_{\max(1,p/n)}(\Hcal)$.

 For all $n\in \N$ and  $x=(x_n,\ldots, x_1)\in F$  we set
\begin{eqnarray}
\label{q-n-x}
q(n;x)&:=&
\sum_{k=1}^n q_{x_k}+\pi
\sum_{k=1}^n \sum_{j=1}^{n-k} \hermsp{\B^{j-1}b_{x_{j+k}}}{a_{x_k}}
,\\
\label{ab-n-x}
 a(n;x)&:=&\sum_{k=1}^n (\B^{n-k})^* a_{x_k} \quad \mbox{ and }\quad
  b(n;x):=\sum_{j=0}^{n-1} \B^{j}b_{x_{j+1}}.
 \end{eqnarray}
For  $n\in \N$ such that $\B^n\in \Scal_2(\Hcal)$ with $\|\B^n\|<1$   we define
(depending on $a_x,\, b_x\in \Hcal$, and $q_x\in \C$) a   function  $c(n;\cdot):F^n\to\R$   via
\begin{equation}
    \label{c-n-x}
c(n;x):=
 \textstyle
\frac{\exp\!\big( 2 \Re{q(n;x)}+ \pi\|a(n;x)\|^2+  \pi \big\|   (1-\B^n (\B^n)^*)^{-1/2}
 \B^n (a(n;x)+b(n;x))
 \big\|^2 \big) }{\det(1-\B^n(\B^n)^*)}.
\end{equation}

 The following proposition describes the type of
operators we will use to build our dynamical
trace formula:

\begin{Prop}
    \label{mixed-iterates-prop}
    Let $F\neq \emptyset$ be an index set,
 $\Hcal$ be a Hilbert space, and $\B\in \End(\Hcal)$.
Given $a_x,\, b_x\in \Hcal$, and $q_x\in \C$ for $x\in F$, consider the operator
$\Mcal_{x}\colon \Ccal(\Hcal)\to \Ccal(\Hcal)$ defined by
\begin{equation}
    \label{mxop-eq}
(\Mcal_{x} f)(z):= \exp\!\big( q_x+ \pi  \hermsp{ z}{a_x} \big)\,
f(b_x+\B z).
\end{equation}
Fix $n\in \N$ and $x_1,\ldots, x_{n}\in F$. Then
$$
(\Mcal_{x_n}\circ \ldots \circ \Mcal_{x_1}f)(z)
=
\exp\!\big(q(n;x)  + \pi \hermsp{ z }{a(n;x)} \big)\,
f\big(\B^n z+b(n;x) \big).
$$
If $\B^n\in \Scal_2(\Hcal)$
with  $\|\B^n\|<1$, then
$\Mcal_{x_n}\circ \ldots \circ \Mcal_{x_1}$ restricts to a Hilbert-Schmidt operator
on $\Fcal(\Hcal)$ which satisfies
$$
\|\Mcal_{x_n}\circ \ldots\circ \Mcal_{x_1} \|^2_{\Scal_2(\Fcal(\Hcal))}\
= c(n;x).
$$
  \end{Prop}
 \begin{proof}
For any family $\{\Mcal_x\mid x\in F \}$ of weighted composition operators
acting via
$(\Mcal_x f)(z)=\exp\!\big( A_x(z)\big)\, f(\psi_x(z))$, induction shows that
$$
(\Mcal_{x_n}\circ \ldots \circ \Mcal_{x_1}f)(z)
=
\exp\Big(
\sum_{k=1}^n (A_{x_k}\circ \psi_{x_{k+1}}\circ\ldots\circ\psi_{x_n})(z)\Big)
\,
(f\circ\psi_{x_1}\circ\ldots\circ \psi_{x_n})(z)
$$
 for all $x_1,\ldots ,x_n\in F$. In particular, if
 $A_x(z)=q_x+\pi \hermsp{z}{a_x}$ and
 $\psi_x(z)=b_x+\B z$ for some $a_x,\, b_x\in \Hcal$, $q_x\in \C$, and $\B\in \End(\Hcal)$
as above, we have to consider   mixed iterates of affine maps:
    Let $V$ be a complex vector space, $\B:V\to V$ a  linear operator, and $b_x\in V$ for $x\in F$.
Then induction shows that
$\psi_{x}:V\to V, \ z\mapsto b_x+\B z$ satisfies
 $$
(\psi_{x_{1}}\circ\ldots\circ\psi_{x_k})(z)=\B^k z+\sum_{j=0}^{k-1} \B^{j}b_{x_{j+1}}
=\B^k z+b(k; (x_k, \ldots, x_1))
$$
for all $k\in \N$, $x_1,\ldots, x_k\in  F$, and $z\in V$.
This implies for $x=(x_n,\ldots, x_1)\in F^n$ that
\begin{eqnarray*}
\lefteqn{
(\Mcal_{x_n}\circ \ldots \circ \Mcal_{x_1}f)(z) \ =
}\\
&=&
\exp\Big(
\sum_{k=1}^n (A_{x_k}\circ \psi_{x_{k+1}}\circ\ldots\circ\psi_{x_n})(z)\Big)
\,
(f\circ\psi_{x_1}\circ\ldots\circ \psi_{x_n})(z)\\
&=&
\exp\Big(
\sum_{k=1}^n A_{x_k} (\B^{n-k} z+\sum_{j=1}^{n-k} \B^{j-1}b_{x_{j+k}}) \Big)
\,
f\big(\B^n z+b(n;x) \big) \\
&=&
\exp\Big(
q(n;x)+\pi
\sum_{k=1}^n \hermsp{\B^{n-k} z}{a_{x_k}}   \Big)
\,
f\big(\B^n z+b(n;x)  \big)
 \\
&=&
\exp\!\big(
q(n;x)   + \pi
 \hermsp{ z }{a(n;x)}\big)\,
f\big(\B^n z+b(n;x)  \big),
\end{eqnarray*}
which yields  the first claim.
Now Corollary~\ref{comp-positive-infty-cor} gives the stated Hilbert-Schmidt norm formula
 and the invariance of  $\Fcal(\Hcal)$.
  \end{proof}

Recall the  matrix subshift $(F,\nu,\A,A) $  with interaction function $A$ from (\ref{a-def})
and the operators $\Mcal_{x}$  defined via
(\ref{mxop-eq}). The following theorem provides a dynamical trace formula for the corresponding
partition function $Z_n(A)$ defined via (\ref{partition-fct}).

  \begin{Theorem}
    \label{thm-1}
Suppose there exists $n_o\in \N$ such that
$
 c(n_o;\cdot )
$ is $\nu^{n_o}$-integrable.
 Then there exists an index $n_1\in \N$ such that
for all $n\geq n_1$ the iterates $\Mcal^n$ of  the Ruelle-Mayer type transfer operator
$$
(\Mcal f)(\sigma, z)
= \int_F \A(x,\sigma)\, \exp\!\big( q_x + \pi \hermsp{z}{a_x}\big)\, f(b_x+\B z)\, d\nu(x)
$$
are trace class operators $ \Mcal^n \in \End \big(L^2(F,d\nu)\hat{\otimes}\Fcal(\Hcal)\big)$.
Moreover, the dynamical partition function $  Z_n(A)$  can be expressed as
$$  Z_n(A)= \det(1-\B^n)\, \trace \Mcal^n    .$$
  \end{Theorem}

\begin{proof}
  By Lemma~\ref{lemma-1} one has for $\tilde{\Hcal}:= L^2(F,d\nu)\hat{\otimes}\Fcal(\Hcal)$
  and $n\geq n_o$
 \begin{eqnarray*}
\|\Mcal^n\|^2_{\Scal_2( \tilde{\Hcal})}
&=& \int_F \int_{F^n}  \A(x_2,x_1)\ldots \A(x_n,x_{n-1})\, \A(\sigma, x_n)\,\times\\
&& \phantom{m}\times\ \| \Mcal_{x_n}\circ \ldots\circ \Mcal_{x_1} \|^2_{\Scal_2(\Fcal(\Hcal))} \,
    d\nu^n(x_1,\ldots ,x_n) \, d\nu(\sigma)\\
&\leq &\nu(F)\,  \int_{F^n}
\|  \Mcal_{x_n}\circ \ldots\circ \Mcal_{x_1} \|^2_{\Scal_2(\Fcal(\Hcal))} \, d\nu^n(x_1,\ldots ,x_n).
  \end{eqnarray*}
There exists $n_1\in \N$ such that for all $n\geq n_1$
  we have $\B^n\in \Scal_1(\Hcal)$ and $\|\B^n\|<1$.
  By Proposition~\ref{mixed-iterates-prop}  (and, possibly,  by enlarging $n_1$)
  the operator $\Mcal^n  $ is trace class for all $n\geq n_1$.
  By Lemma~\ref{lemma-T-HS} the trace of $\Mcal^n$ is given by
   \begin{equation}
   \label{mtilde-trace}
 \trace \Mcal^n   =
\int_{F^n}\! \Big(\prod_{j=1}^{n} \A(x_{j},x_{j+1})\Big)\,
\trace \!(\Mcal_{x_n}\!\circ \ldots\circ \Mcal_{x_1}) \,  d\nu^n(x_1,\ldots, x_n)
   \end{equation}
   using the convention that $x_{n+1}=x_1$.
Proposition~\ref{mixed-iterates-prop} shows that for all choices of
$x=(x_n,\ldots, x_1)\in F^n$ the trace formula from
 Theorem \ref{spurformel-fockinfty-thm} can be applied to the operator
 $\Mcal_{x_n}\circ \ldots\circ \Mcal_{x_1}$. It yields
    \begin{equation}
   \label{mx-trace}
\trace (\Mcal_{x_n}\circ \ldots\circ \Mcal_{x_1}) \
= \
   \frac{\exp\big( \textstyle q(n;x)+ \pi
 \langle (1-\B^n)^{-1} b(n;x)   \mid a(n;x) \rangle  \big)}{\det(1-\B^n)}
    \end{equation}
 with $a(n;x), \, b(n;x)$, and $q(n;x)$ as in  (\ref{q-n-x}) and (\ref{ab-n-x}).
The   inner product occurring in (\ref{mx-trace}) can be rewritten as
\begin{eqnarray*}
 \langle (1-\B^n)^{-1} b(n;x)   \mid a(n;x) \rangle \!\!
& =&
\!\!\big\langle (1-\B^n)^{-1} \sum_{j=0}^{n-1} \B^{j}b_{x_{j+1}}  \mid \sum_{k=1}^n (\B^{n-k})^* a_{x_k} \big\rangle
 \\
 &=& \!\! \sum_{j=0}^{n-1} \sum_{k=1}^n
 \hermsp{ (1-\B^n)^{-1}   \B^{j+n-k}b_{x_{j+1}}  }{  a_{x_k}}
 \\
 &=&\!\!
\sum_{k=1}^n
\sum_{j=0}^{n-1}  \sum_{l=0}^\infty
 \hermsp{   \B^{j+n-k+ln}b_{x_{j+1}}  }{  a_{x_k}} .
   \end{eqnarray*}
Extending  the finite sequence $b_{x_1}, \ldots, b_{x_n}\in \Hcal$ to an
  $n$-periodic sequence, i.\!~e., setting  $b_{x_{k+ln}}:=b_{x_k}$ for all $k=1,\ldots,n$ and $l\in \No$, we obtain
\begin{eqnarray*}
\lefteqn{
q(n;x)+ \pi
 \langle (1-\B^n)^{-1} b(n;x)   \mid a(n;x) \rangle
 =}\\
 &=&\!\! \!
\sum_{k=1}^n q_{x_k}+\pi
\sum_{k=1}^n \sum_{j=1}^{n-k} \langle {\B^{j-1}b_{x_{j+k}}} \mid {a_{x_k}} \rangle
+\pi
\sum_{k=1}^n
\sum_{j=0}^{n-1}  \sum_{l=0}^\infty
 \langle{   \B^{j+n-k+ln}b_{x_{j+1}}  }\mid {  a_{x_k}} \rangle
 \\
 &=& \!\!\!
\sum_{k=1}^n q_{x_k}+\pi
\sum_{k=1}^n
\sum_{j=1}^{\infty} \hermsp{\B^{j-1}b_{x_{j+k}}}{a_{x_k}}
 \\
 &=& \!\!\!
 \sum_{k=0}^{n-1}  A(\tau^k ( \overline{x_1\ldots x_{n}} )) .
    \end{eqnarray*}
This together with (\ref{mtilde-trace}), (\ref{mx-trace}), and (\ref{partition-fct}) completes the proof.
  \end{proof}
We note that unless $\B\in \Scal_2(\Hcal)$ with $\|\B\|<1$, we cannot show
$ \Mcal$ to be a bounded operator on $L^2(F,d\nu)\hat{\otimes}\Fcal(\Hcal)$.

  \section{Ising Type Interactions}
  \label{ising-type-sec}

In this section  we present some explicit models known from the physics literature
for which Theorem \ref{thm-1}
actually provides dynamical trace formulae. A detailed description of these models can be found in
\cite{R06}. Another application of  Theorem \ref{thm-1} are hard rod type models, which are studied in ~\cite{R06} and \cite{R07}.

We consider  a matrix subshift
 $(F,\nu,\A,A) $, where the interaction
  is of the form $\beta\, A_{(q,r,d)}$ with
\begin{equation}\label{aphi-intro}
A_{(q,r,d)}: \Omega_\A \to \C, \ \xi\mapsto q(\xi_1)+\sum_{i=2}^\infty r(\xi_1,\xi_i)\, d(i-1)
.\end{equation}
Here
  $r:F\times F\to\C$ is called an \textit{interaction matrix,} $d\in \ell^1\N$ a \textit{distance function},
  and $q\in \Ccal_b(F)$ a \textit{potential}. The extra parameter $\beta\in \C$ is usually called the
  \textit{inverse temperature}.
 In this context the dynamical  partition function
 $ Z_n(\beta\, A_{(q,r,d)})$
 as defined in (\ref{partition-fct}) coincides with the usual partition function for the
two-body interaction with respect to $(q,r,d)$ and \textit{periodic boundary condition}
(see \cite[Cor 1.11.3]{R06}).
We will assume additional properties of the distance function and the interaction matrix. A list of examples will be given in
Examples~\ref{ex-1} and \ref{ex-2}.

An interaction matrix $r$ is said to be of  \textit{Ising type}, if
  there exists a finite number of  functions $s_i, t_i:F\rightarrow \C$ such that
$$
r(x,y)=\sum_{i=1}^M \overline{s_i(x)} \,t_i(y).
$$
 The minimal number $M$ is called the \textit{rank} of $r$.

\begin{Example}
    \label{ex-1}
\begin{enumerate}
\item
\textit{Ising model:} Let $F\subset \C$ be a bounded set and
$\ds r(x,y)=xy.
$
In E. Ising's  original model \cite{I25} he took  $F=\{\pm 1\}$, the so called
\textit{spin-$\einhalb$ model}, in order to describe
ferromagnetism of a solid, where the spins of the electrons can only take values in a set with two elements,
``spin up'' or ``spin down''.

\item
Let $F\subset \Hcal$ be a bounded subset of a finite dimensional Hilbert space $(\Hcal, \hermsp{\cdot}{\cdot})$. Then
any bilinear form $\beta$ on $\Hcal$ defines an interaction matrix of Ising type via
$
r(x,y):= \beta(x,y)  $. In fact: Choose an orthonormal basis $(e_j)_{j=1, \ldots, \dim \Hcal}$ for $\Hcal$, then
$$  \beta(x,y)
=\beta(\sum_{i=1}^{\dim \Hcal}  \hermsp{x}{e_i}e_i, y)
=\sum_{i=1}^{\dim \Hcal}    \hermsp{x}{e_i}  \beta(e_i,y)
.
$$
  Note that  $r$ has  rank less or equal to $\dim \Hcal$.

\item
The (generalized)
\textit{Stanley $M$-vector model,} cf. \cite{St68},
is a special case of (ii):  Take  $\Hcal=\R^M$ with $r(x,y)=  \hermsp{x}{y}
  $  and
$
F:=
 \big\{ \sigma \in \R^M\,\mid \ r(\sigma,\sigma)=s^2\big\}$  to be the   $(M-1)$-sphere   with radius $s>0$
equipped with the (normalized) surface measure $\nu$ on $F$.

The following  table gives a list of physical models which can be seen as
 applications of Stanley's $M$-vector model. Depending on the parameter
$M$ these models have special names.
\end{enumerate}
\begin{flushright}
\begin{tabular}[t]{|c|l|l|}\hline
Rank& Special name & System \\ \hline\hline
$1$& Ising model&  one-component fluid, binary alloy, mixture\\ \hline
$2$& Planar model
&  $\lambda$-transition in a Bose fluid\\ \hline
$3$& Heisenberg model&  (anti-)ferromagnetism\\ \hline
$M>3$&   $M$-vector  model& no physical system discovered yet\\ \hline
\end{tabular}
\end{flushright}
\begin{enumerate}
\item[]
This table is taken from \cite[p.~488]{St74} where one can also find a lot of references to the underlying
physics.
Note that the rank $1$ case gives $F=\{\pm 1\}$ and hence the spin-$\einhalb$ Ising model.
\item[(iv)]
   If $F$ is finite, then \textit{every} interaction matrix is of Ising type, since
$$
r(x,y)=\sum_{z\in F}  r(x,z)\, \delta(y,z),
$$
where $\delta:F\times F\rightarrow \C$ is  Kronecker's delta on $F$.
In particular, the finite-state \textit{Potts model} is of Ising type:
Let $F$ be a finite set and
$\ds r(x,y)=\delta(x,y) $.
This
  model is due to R. Potts~\cite{Po52} and describes the situation where only electrons with identical spin
interact.
 \qed
\end{enumerate}
\end{Example}
We remark that Ising, Potts, and Stanley have considered these models only for finite range interactions.

Consider the interactions $A_{(q, r,d )}$ with distance functions $d$ belonging to
 subspaces $\Dcal_p \subset \ell^1\N$ (for $p\in[1,\infty[$) which are defined as follows:
$d=d_{\B,v,w}\in \Dcal_p $ if and only if there exist  a Hilbert space $\Hcal$ and a
Schatten class operator
 $\B \in  \Scal_p(\Hcal)$ with  spectral radius $\rhospec(\B)<1$ and vectors
  $v,\,w\in \Hcal$ such that
$$ d:\N\to\C,\ k\mapsto d(k)=\hermsp{ \B^{k-1}v}{w}_{\Hcal}.
$$

Now, let $d=d_{\B,v,w}\in \Dcal_p$ and $r(x,y)=\sum_{i=1}^M \overline{s_i(x)} \,t_i(y)
$. We rewrite the  Ising type observable
$
A_{(q,r,d)}: \Omega_\A\to \C$
in such a way that we can apply  Theorem~\ref{thm-1}.
\begin{eqnarray*}
A_{(q,r,d)}(\xi)
&=&
 q(\xi_1)+\sum_{i=2}^\infty \sum_{j=1}^M \overline{s_j(\xi_1)}\,t_j(\xi_i)\, \hermsp{  \B^{i-2}v }{ w }_{\Hcal}
  \\
 & =&
 q(\xi_1)+\sum_{i=2}^\infty \hermsp{ (\underline{\B_M})^{i-2} \underline{t_M}(\xi_i,v)  }{   \underline{s_M}( \xi_1,w)      }_{\Hcal^M}
 ,
\end{eqnarray*}
where $\underline{\B_M}:\Hcal^M   \to \Hcal^M,\, (z_1,\ldots, z_M)\mapsto (\B z_1,\ldots, \B z_M)$,
$\underline{s_M}:F\times \Hcal \to \Hcal^M, \, \underline{s_M}(x,v):=  \big(s_1(x)v,\ldots, s_M(x)v\big), $ and, similarly,  $\underline{t_M}:F\times \Hcal \to \Hcal^M,$ $  \underline{t_M}(x,v):=  \big( t_1(x)v,\ldots, t_M(x)v\big)$.

\begin{Theorem}
     \label{thm-2}
    Let
 $(F,\nu,\A,A_{(q,r,d)}) $  be a matrix subshift  where
  $q\in\Ccal_b(F)$,
    $d=d_{\B,v,w}\in \Dcal_p$ and
      $r\in \Ccal_b(F\times F)$ is an interaction matrix of Ising type, say
$
r(x,y)=\sum_{i=1}^M \overline{s_i(x)}\, t_i(y)$
with $s_i,\, t_j\in \Ccal_b(F)$.
Then there exists an index  $n_0\in\N$ depending on $\B$ such that for all $n\geq n_0$ the iterates
$
\Mcal_\beta^n\in \End( L^2(F,d\nu)\hat{\otimes}\Fcal(\Hcal^M) )$
of the
Ruelle-Mayer transfer operator
$$
(\Mcal_\beta f)(x,z )=
\int_F \A_{\sigma,x}\exp\!\Big(\beta q(\sigma)+ \beta \langle z| \underline{s_M}(\sigma ,w)\rangle\Big)\,
f(\sigma, \underline{t_M}(\sigma,v) +\underline{\B_M} z)
 \, d\nu(\sigma)
$$
are of trace class and satisfy the dynamical trace formula
$$
\zn(\beta A_{(q,r,d)}) =\det(1-\B^n)^M\,
\trace \Mcal_\beta^n .$$
\end{Theorem}

\begin{proof}
 By assumption the sets
 $\{a_x:=\underline{s_M}(x,w)\in \Hcal^M\,|\, x\in F\} $ and
 $\{b_x:=\underline{t_M}(x,w)\in \Hcal^M\,|\, x\in F\} $  are bounded.  Choose $m\in \N$
 large enough such that
 $\B^{m}$  has operator norm  less than one and is Hilbert-Schmidt.
 Proposition~\ref{mixed-iterates-prop}, applied to $\underline{\B_M}$, shows that
  the associated function $c(m;\cdot)$ is bounded, hence integrable. Thus we can
 apply  Theorem~\ref{thm-1} to $\Mcal_\beta$  which  proves the claim.
\end{proof}

\begin{Corollary}[Ising model]
    \label{thm-1-cor2}
    Let $F\subset \C$ be a bounded set equipped with a finite measure $\nu$
and
 $(F,\nu,\A,A_{(q,r,d)}) $  be a matrix subshift,  where
  $q\in\Ccal_b(F)$,
    $d=d_{\B,v,w}\in \Dcal_p$ and $r(x,y)=xy$.
Then there exists an index  $n_0\in\N$ depending on $\B$ such that for all $n\geq n_0$ the iterates
$
\Mcal_\beta^n\in \End( L^2(F,d\nu)\hat{\otimes}\Fcal(\Hcal) )$
of the
Ruelle-Mayer transfer operator
$$
(\Mcal_\beta f)(x,z)=
\int_F \A_{\sigma,x}\, \exp\!\Big(\beta q(\sigma)+ \beta \sigma
\langle z|w\rangle\Big)\,
f(\sigma, \sigma\,v+\B\,z)\, d\nu(\sigma)
$$
are trace class and
satisfy
$
\zn(\beta A_{(q,r,d)})
=\det(1-\B^n)\,
\trace \Mcal_\beta^n.
$
\qed
\end{Corollary}
\begin{Corollary}[Potts model]
  Let $F=\{1,\ldots,N\}$
be a finite alphabet,  the measure $\nu$ on $F$ be identified with  its
distribution vector, and  $(F,\nu,\A,A_{(q,r,d)}) $  be a matrix subshift,  where
   $q:F\to\C$,
    $d=d_{\B,v,w}\in \Dcal_p$ and $r(x,y)=\delta(x,y)$.
Then there exists an index  $n_0\in\N$ depending on $\B$ such that for all $n\geq n_0$ the iterates
$
\Mcal_\beta^n\in \End( L^2(F,d\nu)\hat{\otimes}\Fcal(\Hcal^N) )$
of the
Ruelle-Mayer transfer operator
\begin{eqnarray*}
\lefteqn{
(\Mcal_\beta f)(l; z_1,\ldots,z_N)\, =
}
\\&&
=\sum_{k=1}^N \A_{k,l}\,\nu_k\,
 \exp\!\Big(\beta\, q(k)+ \beta\, \langle z_k|w\rangle\Big)\,
f\big(k; ( \delta_{k,m}v +\B z_m)_{m=1,\ldots,N}\big).
\end{eqnarray*}
are trace class and
satisfy
$
\zn(\beta A_{(q,r,d)})
=\det(1-\B^n)^N\,
\trace \Mcal_\beta^n.
$
\qed
\end{Corollary}

By the canonical identification       $ L^2(F,d\nu)\hat{\otimes}\Fcal(\Hcal^N) \cong \Fcal(\Hcal^N)^{N}$
the Ruelle-Mayer transfer operator can be rewritten as
\begin{eqnarray*}
\lefteqn{
((\Mcal_\beta (f_1,\ldots, f_N))_l(z_1,\ldots,z_N)\, =
}
\\&&
=\sum_{k=1}^N \A_{k,l}\,\nu_k\,
 \exp\!\Big(\beta\, q(k)+ \beta\, \langle z_k|w\rangle\Big)\,
f_k\big(( \delta_{k,m}v +\B z_m)_{m=1,\ldots,N}\big).
\end{eqnarray*}

\begin{Example}
    \label{ex-2}
    \label{classification-example}%
\begin{enumerate}
\item
Finite range: There exists $\rho_0\in\N$, the range of $d$,  such that
$\ds \, d(k)=0 \mbox{ for all } k>\rho_0$. Remark~\ref{srho-prop} below shows that $d\in \Dcal_1$.

\item
Polynomial-exponential:
$\ds\,
d:\N\to\C,\, k\mapsto \lambda^k\, p(k),$
where $p\in \C[z]$ is a polynomial and $\lambda\in \C$ with  $0<|\lambda|<1$ is the decay rate.
Remark~\ref{b-prop} below shows that $d\in \Dcal_1$.

\item
Superexponential:
Let
 $\gamma>0,\, \delta>1$ and
$\, \ds
 d:\N\to\C,\, k\mapsto
a(k)\,  \exp(-\gamma\,k^{\delta}),
$
where
$a:\N\to\C$
 is  of lower order
 such  that   $\ds\lim_{k\to\infty} a(k)\, \exp(-\epsilon_1k^{\epsilon_2})=0
 $ for all $\epsilon_1,\,\epsilon_2>0$.
Proposition~\ref{sg-prop} below shows that $d\in \Dcal_1$. (The decay estimate can be weakened, cf. Example~\ref{super-exp-example}.)

\item
Suitable infinite superpositions of exponentially decaying terms:
$$
d(k)=\sum_{i=1}^\infty c_i\, \lambda_i^k,$$
where $\lambda\in \ell^p\N$ ($1\leq p <\infty$) and  $c:\N\to\C$
such that $c\lambda:\N\to \C,\, n\mapsto c_n\,\lambda_n$ belongs to $ \ell^1\N$.
Obviously, $d\in \Dcal_p$.

Example: Pick a holomorphic function $f$  on the unit disk with $f(0)=0$, $0<|\lambda|<1$, then $d(k)=f(\lambda^k)$  belongs to $\Dcal_1$.
\qed
\end{enumerate}
\end{Example}

We conclude this section with the verification of the various claims made in
Example \ref{classification-example}.
We start with  some obvious facts on  finite-range interactions.

\begin{Remark}
    \label{srho-prop}%
For $\rho_0\in \N_{>1}$ let
\begin{equation*}
 \S_{\rho_0}:=
\left(\!\!\begin{array}{cccc}
0&1\\
&\ddots&\ddots\\
&&0&1\\
&&&0
\end{array}\!\!\right)\in \Mat({\rho_0},{\rho_0};\Z).
\end{equation*}
\begin{enumerate}
\item
Then $\S_{\rho_0}$
is a $\rho_0$-step nilpotent matrix with spectral radius $\rhospec(\S_{\rho_0})=0$.
\item
Let
 $d:\N\rightarrow \C$ be a finite range distance function, say
 $d(k)=0$ for all $k>{\rho_0}$, $\lambda\in \C^\times$, and  $w^d\in\C^{\rho_0}$ with entries
 $w^d(k)=\lambda^{1-k}\,d(k)$. Then
$
d(k)
= \hermsp{  (\lambda\S_{\rho_0})^{k-1}w^d }{ e_1 }
$
for all $k\in \N$, where $e_1=(1,0,\ldots,0)$.
\qed
\end{enumerate}

\end{Remark}

Fix an additional decay parameter $\lambda$ with $|\lambda|<1$. Then
$\|( \lambda \S_{\rho_0})^k\|<1$ for all $k\in \N$ and the
 dynamical trace formula for the Ruelle-Mayer transfer operator built from  $\lambda \S_{\rho_0}$
 holds for all $n\in \N$ (instead just for  $n\geq \rho_0$ in the case $\lambda=1$).

\begin{Remark}
    \label{b-prop}
Let  $p\in \No$ and  define $v:\R\to \R^{p+1}$ via
$v(x):=\lambda \big( 1, (1+x), \ldots, (1+x)^p\big)$.
Further, let $\B^{(p+1)}\in \Mat(p+1,p+1;\R)$ be the unipotent (lower) triangular matrix with entries
 $$(\B^{(p+1)})_{i,j}=\left\{\begin{array}{cl}
{ i\choose j} &, j\leq i,\\
0&, \mbox{otherwise.}
\end{array}\right.
$$
Then $\det(1-(\lambda \B^{(p+1)})^n)=(1-\lambda^n)^{p+1}$ and the binomial formula implies
$$
\left(\!\!\begin{array}{c}
\lambda^{k+1}\\
(k+1+x)\lambda^{k+1}\\
\vdots\\(k+1+x)^{p-1}\lambda^{k+1}\\
(k+1+x)^p\lambda^{k+1}
\end{array}\!\!\right)=\lambda
\B^{(p+1)}
\left(\!\!\begin{array}{c}
\lambda^{k}\\
(k+x)\lambda^{k}\\
\vdots
\\
(k+x)^{p-1}\lambda^{k}\\
(k+x)^p\lambda^{k}
\end{array}\!\!\right)
$$
for all $k\in \N$, and induction yields
$$
\left(\!\!\begin{array}{c}
\lambda^{k}\\
(k+x)\lambda^{k}\\
\vdots
\\
(k+x)^{p-1}\lambda^{k}\\
(k+x)^p\lambda^{k}
\end{array}\!\!\right) =
(\lambda \B^{(p+1)})^{k-1} v(x).
$$
Consequently, for $c:= \big( c_0,  \ldots, c_p\big) \in \C^{p+1}$ and $x=0$ we have
$$
\lambda^k\,\sum_{i=0}^p c_i\, k^i=\hermsp{ \lambda^k
(\B^{(p+1)})^{k-1} \underline{1} }{ \overline{c} }
$$
with $\underline{1}=(1,\ldots,1)\in \Z^{p+1}.$
Thus the  distance functions
$d(k):=\lambda^k\,\sum_{i=0}^p c_i\, k^i$ belongs to $\Dcal_1$.
\qed
\end{Remark}

\begin{Prop}
    \label{sg-prop}%
Let $g:\N\to \C\setminus\{0\}$ with  $\sum_{k=1}^\infty
\big|\frac{g(k)}{g(k+1)}\big|^p<\infty.$
We define
$\S_g:\C^\N\to \C^\N,\ (\S_g z)_k
:=\frac{g(k)}{g(k+1)}\, z_{k+1}$.
 Then:
\begin{enumerate}
\item
$\S_g$ leaves  the spaces $\ell^q\N$ invariant
for $1\leq q<\infty$. Moreover, it defines continuous operators  with
$
\|\S_g \|_{\ell^q\N\to \ell^q\N}\leq \sup_{k\in \N}
\big|\frac{g(k)}{g(k+1)}\big|
$
on these spaces.

\item
 For all $z\in \ell^p\N$, $n\in\No$, $k\in \N$, we have
$
(\S_g^{n} z)_k=
\frac{g(k)}{g(k+n)}\, z_{k+n}.
$

\item
 $\S_g:\ell^2\N\to\ell^2\N $ belongs to the Schatten  class $\Scal_p(\ell^2\N)$.
 It is not normal, and it satisfies    $\rhospec(\S_g)<1$ as well as
 $\det(1-\S_g^n)=1$ for all $n\geq p$.

\item
If, for  $d:\N\to\C$, the function
$v^d:\N\to\C,\ v_k^d:=g(k) \, d(k)$ belongs to $ \ell^2\N$, then
$$
d(k)=\frac{1}{g(1)}\,
\hermsp{ \S_g^{k-1} v^d }{ e_1 }_{\ell^2\N},
$$
where $e_1=(1,0,\ldots,0)$.
 In particular $d=d_{\S_g,v^d,e_1}\in  \Dcal_p$ for $1\le p <\infty$.
 \end{enumerate}
\end{Prop}

\begin{proof}
Let $1\leq q<\infty$ and  $z\in \ell^q\N$, then
$$
\|\S_g z\|_{\ell^q\N}^q\,=\,
\sum_{k=1}^{\infty}
\Big|\frac{g(k)}{g(k+1)}\, z_{k+1}\Big|^q
\,\leq \, \sup_{k\in \N}
\Big|\frac{g(k)}{g(k+1)}\Big|^q\, \|z\|_{\ell^q\N}^q.
$$
This implies (i).
Assertion (ii) is easily shown by induction.
The
$\ell^2\N$-adjoint $\S^*_g$ of $\S_g$ is given by
$$
\S_g^*:\ell^2\N\to \ell^2\N,\
 (\S_g^* \xi)_i=\left\{\begin{array}{ll}0,&i=1,\\
\overline{
\frac{g(i-1)}{g(i)}}\, \xi_{i-1}, &i\geq 2.\end{array}\right.
$$
Therefore
$
((\S_g^{}\S_g^*)(\xi))_k=
\big|\frac{g(k)}{g(k+1)}\big|^2\, \xi_k,
$
which shows that $\S_g^{}\S_g^*$ is diagonal with respect to the standard basis.
We can read off the singular numbers of $\S_g$ being the  square roots of the diagonal entries of
$\S_g^{}\S_g^*$. By assumption they belong to $\ell^p\N$.
On the other hand
$$
((\S^*_g\S_g^{})(\xi))_k=\left\{\begin{array}{ll}0,&i=1,\\
\Big|\frac{g(i-1)}{g(i)}\Big|^2\, \xi_{i}, &i\geq 2.\end{array}\right.
$$
The operator norm of $\S^n_g$ is bounded by
$\sup_{k\in \N}\big|\frac{g(k)}{g(k+n)}\big|$.
The sequence $k\mapsto \big|\frac{g(k)}{g(k+1)}\big|$ tends to zero, hence one can find $k_0\in \N$
such that
$ \big|\frac{g(k)}{g(k+1)}\big|<\einhalb$ for all $k\geq k_0$.
    Let $C=\max_{k=1,\ldots, k_0} \big|\frac{g(k)}{g(k+1)}\big|$.
Then for all $k\in \N$ one has
\begin{eqnarray*}
 \Big|\frac{g(k)}{g(k+n)}\Big|
 &=&
  \Big|\frac{g(k)}{g(k+1)}\frac{g(k+1)}{g(k+2)}\cdots
  \frac{g(k_0)}{g(k_0+1)}\frac{g(k_0+1)}{g(k_0+2)}\cdots
  \frac{g(k+n-1)}{g(k+n)}\Big|
  \\&\leq&
2^{-(n-k_0+k)} \, C^{k_0-1},
  \end{eqnarray*}
  which tends to zero as $n\to \infty$.
Hence we can find $n\in \N$ such that $\|\S^n_g\|<1$ and hence
$\rhospec(\S_g)<1$.
With respect to the standard basis of $\ell^2\N$ the operator $\S_g$ is an upper triangular matrix with
 zeros along the diagonal, hence
 $\det(1-\S^n_g)=1$ for all $n\geq n_0$.
As a consequence of   (ii) we have
$$
(\S^n_g v^d)_l
=\frac{g(k)}{g(l+n)}\,v^d_{l+n}
=\frac{g(l)}{g(l+n)}\,g(l+n)\, d(l+n)
$$
for all $n\in\No$, $l\in\N$, which immediately  implies that
$$
\hermsp{ \S^{k-1}_g v^d }{ e_1 }=
(\S^{k-1}_g v^d)_1
=g(1)\, d(k).
$$
\end{proof}

\begin{Remark}
    \label{super-exp-example}%
\begin{enumerate}
\item
For any nowhere-vanishing sequence $s\in\ell^p\N$ setting
$g:\N\to\C, \ g(k):=\big(\prod_{l=1}^{k-1}s(l)\big)^{-1}$ one obtains
a function $g$ of the kind required in Proposition~\ref{sg-prop}. In particular, $
s(k)=\frac{g(k)}{g(k+1)}$ and $|s|:\N\to\C,\ n\mapsto |s(n)|$
is the sequence of singular numbers of the corresponding weighted shift operator $\S_g$.
Functions $g:\N\to\C$ satisfying  the summability condition
 $\sum_{k=1}^\infty \big|\frac{g(k)}{g(k+1)}\big|^p<\infty$
are, for instance,
$g(k):=\exp(\gamma k^\delta)$ with $\gamma>0,\, \delta>1$.

\item
More generally, consider the following distance function
$d:\N\to\C$ which was studied by   D. Mayer in  \cite[p. 100]{M80a}:
 $d(k)= a(k)\,  \exp(-\gamma\,k^{\delta}),$
where
 $\gamma>0,\, \delta>1$ and $a:\N\to\C$
 is  a lower order term, in the sense  that
 $\lim_{k\to\infty} a(k)\, \exp(-\epsilon_1\,k^{\epsilon_2})=0
 $ for all $\epsilon_1,\,\epsilon_2>0$ .
We claim that
$$
d(k)=
\hermsp{ \S_g^{k-1} v^d }{ e_1 }_{\ell^2\N},
$$
where
$v^d(k):=a(k)\,  \exp\big(\gamma ((k-1)^\delta-k^\delta)\big) $
 defines $v^d\in \ell^2\N$, and
  $g:\N\to \C,\ k\mapsto \exp(\gamma (k-1)^\delta)$
  defines
  $\S_g \in \Scal_1(\ell^2\N)$ with $\|\S_g\|<1$.

To prove the claim note first that
$g(1)=1$ and $g$ satisfies the summability condition from Proposition~\ref{sg-prop}:
For $\delta>1$ and $j\geq 0$ we have
$$j^\delta-
(1+j)^\delta
=j\,j^{\delta-1}-(1+j)\,(1+j)^{\delta-1}\leq (j-1-j)\, j^{\delta-1}=- j^{\delta-1}.$$
Hence
\begin{equation}
    \label{gg-estimate}%
\sum_{k=1}^\infty
\Big|\frac{\exp(\gamma (k-1)^\delta)}{\exp(\gamma \,k^\delta)}\Big|^p
\leq
 \sum_{ k=1}^\infty
 \exp\!\big(-\gamma\, p\, (k-1)^{\delta-1}\big),
\end{equation}
which is finite
 for all $p>0$.
Hence the corresponding weighted shift operator $\S_g:\ell^2\N\to\ell^2\N $ is trace class.
Moreover,  $\S_g$ has operator norm bounded by $\exp(-\gamma)<1$.
It remains to show that $v^d\in \ell^2\N$. We proceed similar to the previous estimate
(\ref{gg-estimate}).
For $0<\epsilon_1<\gamma $, $0<\epsilon_2\leq \delta-1$,
by our assumptions on  the lower order term $a$  we can find a constant $C>0$ such that
\begin{eqnarray*}
\|v^d\|_{\ell^1\N}
&=&
\sum_{k=1}^\infty
\exp\!\big(-\gamma   (k^\delta-(k-1)^\delta)\big)\,|a(k)|
\\&\leq &
C\,
\sum_{k=1}^\infty
\exp(-\gamma   k^{\delta-1}+\epsilon_1k^{\epsilon_2})
\\&\leq& C\,
\sum_{k=1}^\infty
\exp\!\big(-(\gamma- \epsilon_1)  k^{\delta-1}\big) <\infty.
\end{eqnarray*}
This proves the claim and hence shows that Proposition~\ref{sg-prop}~(iv) applies to
$d$.

 \item
For $d\in \Dcal_q$ for some $q$  it would be sufficient that
$v^d\in \ell^2\N$ and (\ref{gg-estimate}) holds for some $p<\infty$. These observations
allow to weaken the conditions on the lower order term.
For instance,  the sequence $a$ might grow like $k\mapsto \exp(\gamma\,  k^{\delta-1-\epsilon})$
for all $\epsilon >0$.
\qed
\end{enumerate}
\end{Remark}


\end{document}